\documentclass{amsart}
\usepackage[all]{xy}
\usepackage[dvips]{graphicx}
\usepackage{latexsym,amsthm,amssymb,amscd,stmaryrd}
\usepackage[vcentermath]{youngtab}
\usepackage{enumerate}
\usepackage{multicol}
\usepackage{amscd}
\usepackage{amsmath}

\usepackage{floatflt}
\usepackage{color}

\numberwithin{equation}{section}

\newtheorem{theorem}{Theorem}[section]
\newtheorem{lemma}[theorem]{Lemma}

\newtheorem{prop}[theorem]{Proposition}

\newtheorem{ex}[theorem]{Examples}
\newtheorem{exs}[theorem]{Examples}
\newtheorem{conjecture}[theorem]{Conjecture}

\newcommand{\cA}{{\mathcal A}}

\newcommand{\cC}{{\mathcal C}}
\newcommand{\cD}{{\mathcal D}}

\newcommand{\cF}{{\mathcal F}}

\newcommand{\cN}{{\mathcal N}}
\newcommand{\cO}{{\mathcal O}}
\newcommand{\cP}{{\mathcal P}}

\newcommand{\cT}{{\mathcal T}}
\newcommand{\cU}{{\mathcal U}}

\newcommand{\mh}{\mathfrak{h}}
\newcommand{\mb}{\mathfrak{b}}

\newcommand{\mn}{\mathfrak{n}}

\newcommand{\mV}{\mathbb{V}}

\newcommand{\mC}{\mathbb{C}}

\newcommand{\mP}{\mathbb{P}}
\newcommand{\mZ}{\mathbb{Z}}

\newcommand{\la}{\lambda}
\newcommand{\oneone}{\stackrel{1:1}\leftrightarrow}

\newcommand{\HOM}{\operatorname{Hom}}

\newcommand{\END}{\operatorname{End}}

\newcommand{\p}{\mathfrak{p}}

\newcommand{\gMOD}{\operatorname{gmod}}

\begin{document}

\title[A combinatorial approach]{A combinatorial approach to functorial quantum $\mathfrak{sl}_k$ knot invariants}

\author{Volodymyr Mazorchuk}
\address{V. Mazorchuk: Department of Mathematics, Uppsala University (Sweden).}
\email{mazor\symbol{64}math.uu.se}
\author{Catharina Stroppel}
\address{C. Stroppel: Department of Mathematics, University of Glasgow
(UK).} \email{c.stroppel\symbol{64}maths.gla.ac.uk}
\thanks{ V. M. was partly supported by STINT,
the Royal Swedish Academy of Sciences, and the Swedish Research Council.\\
C.S. was partly supported by an EPSRC and a Von-Neumann Award of the
Instutute of Advanced Studies in Princeton. }

\begin{abstract}

This paper contains a categorification of the $\mathfrak{sl}(k)$ link invariant
using parabolic singular blocks of category $\cO$. Our approach is intended to
be as elementary as possible, providing combinatorial proofs of the main
results of \cite{Sussan}. We first construct an exact functor valued invariant
of webs or ``special'' trivalent graphs labelled with $1, 2, k-1, k$ satisfying
the MOY relations. Afterwards we extend it to the $\mathfrak{sl}(k)$-invariant
of links by passing to the derived categories. The approach of \cite{Ksl3}
using foams appears naturally in this context. More generally, we expect that
our approach provides a representation theoretic interpretation of the
$\mathfrak{sl}(k)$-homology, based on foams and the Kapustin-Lie formula, from
\cite{MackaayStosic}. Conjecturally this implies that the Khovanov-Rozansky
link homology is obtained from our invariant by restriction.
\end{abstract}

\maketitle

\date{}
\tableofcontents

\section{Introduction}
Let $k\geq 2$ be a positive integer. In \cite{MOY}, Murakami, Ohtsuki and
Yamada developed a graphical calculus for the $\mathfrak{sl}(k)$ polynomial
invariant ${\mathbf P}_k$ of knots and links. Web diagrams describe
intertwiners between the finite tensor products of fundamental representations
of $\cU=\cU_q(\mathfrak{sl}_k)$, the (generic) quantised universal enveloping
algebra of $\mathfrak{sl}_k$. The $\mathfrak{sl}(k)$ link polynomial ${\bf
P}_k$ is defined via the skein relation
\begin{eqnarray*}
q^k{\mathbf
P}_k\left(\;\parbox[thb]{0.4cm}{\includegraphics{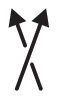}}\;\right)
-q^{-k}{\mathbf
P}_k\left(\;\parbox[thb]{0.4cm}{\includegraphics{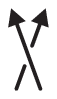}}\;\right)
=(q-q^{-1}){\mathbf{P}_k}\left(\;\parbox[thb]{0.4cm}{
\includegraphics{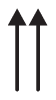}}\;\right)
\end{eqnarray*}
and normalised by setting $\mathbf{P}_k$ of the trivial knot equal to the
quantum number $[k]$.

In this paper we want to describe a categorification of this invariant
${\mathbf P}_k$ using parabolic categories $\cO$ for various $\mathfrak{gl}_n$.
For the special case of $k=3$ we explicitly describe how the
$\mathfrak{sl}_3$-link homology from \cite{Ksl3} emerges naturally from our
approach. More generally, our results should be the representation theoretic
explanation of \cite{MackaayStosic}, which uses foams and the Kapustin-Lie
formula (see Conjecture~\ref{conjecture}). Having set up the representation
theoretic picture conveniently, the verification of this claim reduces to
straight forward, but apparently quite lengthy, combinatorics. In the present
paper, we therefore want to focus on giving all the necessary representation
theoretic tools. Since the Mackaay-Stosic-Vaz homology is equivalent (see
\cite{MackaayStosic}) to the Khovanov-Rozansky homology \cite{KhovRozMatrix},
the verification of the conjecture would give a representation theoretic
interpretation of \cite{KhovRozMatrix}.

In connection with categorifications of link polynomials, in particular the
MOY-relations, category $\cO$ appeared already in several disguises in the
literature. Our results here are a generalisation of \cite{StDuke}, where the
case of the Jones polynomial, i.e. $k=2$, was established. A categorification
for general $k$ using certain (derived categories of) singular blocks of
category $\cO$ was first worked out by Josh Sussan in the paper \cite{Sussan},
which motivated our work. Our picture here will be Koszul dual to Sussan's
(\cite{MOS}). Although very similar on the first sight, our approach appears to
us as being much simpler and better adapted, for instance because of the
following:

\begin{itemize}
  \item The categorification of webs which appears when completely flattening any link
  diagram can be done by working inside certain {\it abelian} categories. Only
  crossings force us to pass to derived categories
  (whereas the approach of \cite{Sussan} has to use derived categories and
  higher derived functors from the very beginning).
  \item Assuming a few standard facts on projective functors turns the problem
  of checking the MOY relations into an easy task, involving a couple of simple facts from the Kazhdan-Lusztig
  combinatorics.
  \item Our approach directly shows the connection to \cite{Ksl3} and
  \cite{MackaayStosic}. The homology rings of partial flag varieties here arise
  as endomorphism rings of projective modules in our picture (using a very
  special and easy case of Soergel's endomorphism theorem \cite{Sperv}).
\end{itemize}

\subsubsection*{The organisation of the paper and the main results}
The main goal
of this paper is to provide a ``down-to-earth'' approach to the quite involved,
technical work of \cite{Sussan}. The price to pay is that one has to assume a
few standard facts about projective functors which we state as {\bf Fact 1} to
{\bf Fact 4} in Section~\ref{Section4}. The MOY-relations are then easy to
check: We first do some calculations in the Hecke algebra of the symmetric
group $S_n$ which describes the combinatorics of projective functors for the
ordinary (non-parabolic) category $\cO$. As a consequence we get the MOY
relations up to some ``error term''.
This ``error term'' vanishes however when we
restrict the functors to the parabolic categories which are really used in our
categorification. Again, the verification is completely combinatorial using the
knowledge of annihilators of induced modules for the symmetric group ({\bf Fact
3}). In fact, only the verification of Reidemeister I
 and one additional move (Proposition~\ref{Sussangap}) involving crossings, require non-combinatorial
 arguments. (Note that the arguments in \cite{Sussan} for these moves
 are incomplete.)

Let now $V$ be the natural representation of $\cU$, i.e. of quantum
$\mathfrak{sl}_k$, and let $\nu$ be a composition of $n$. Consider a tensor
product of fundamental representations of $\cU$ of the form
$$X^\nu:=\bigwedge^{\nu_1}{V}\otimes\bigwedge^{\nu_2}{V}\otimes\ldots\otimes\bigwedge^{\nu_k}{V}.$$
In Section~\ref{Section2} we categorify this $\mC[q,q^{-1}]$-module via the
direct sum
\begin{eqnarray*}
\cC_{X^\nu}:=\bigoplus_\mu {}^\mZ\cO(n)_\nu^\mu
\end{eqnarray*}
of parabolic singular blocks of (the graded version of) category $\cO$ for
$\mathfrak{gl}(n)$, where $\mu$ runs through all compositions of $n$ with at
most $k$ parts. This is a generalisation of the categorifications in
\cite{BFK}, \cite{StDuke}, see also \cite{Brundan}. In
Subsection~\ref{categoryO} we give an explicit isomorphism $\Gamma^\nu$ between
the standard basis vectors of $X^\nu$ and the isomorphism classes of parabolic
Verma modules using some easy combinatorics. This is used afterwards in
Section~\ref{Section4} to categorify intertwiners via graded translation
functors. In Section~\ref{Section4} we show that these translation functors
satisfy the MOY relations for trivalent graphs. This means that to each
``special intertwiner'' $f$ (see Section~\ref{Section2}) labelled by numbers
from $\{1,2,k-1,k\}$ only, we associate in Section~\ref{Section5} some functor
$F(f)=F_k(f)$ such that the following holds:

\begin{theorem}
\label{main1} Let $k\geq 2$ as above and let $\nu$, $\nu'$ be compositions of
$n$.
\begin{enumerate}
\item If $f:X^\nu\rightarrow X^{\nu'}$ is a composition of special intertwiners
then $F(f)$ is an exact functor $\cC_{X^\nu}\rightarrow\cC_{X^{\nu'}}$.
\item Up to isomorphism, the functors satisfy the MOY relations (Figures~\ref{fig:Rel1} to
~\ref{fig:Rel5}).
\item Under the isomorphism $\Gamma^\nu$, a composition $f:X^\nu\rightarrow X^{\nu'}$ of special intertwiners corresponds to $[F(f)]$, the $\mC[q,q^{-1}]$-linear map from the complexified
Grothendieck group of $\cC_{X^\nu}$ to the one of $\cC_{X^{\nu'}}$.
\end{enumerate}
\end{theorem}

In Section~\ref{Section5} we extend this assignment $f\mapsto F(f)$ to a
categorification of the MOY-tangle invariant, by associating to each oriented
tangle diagram $t$ a certain functor $F(t)=F_k(t)$ such that the following
holds:

\begin{theorem}
\label{main2}
\begin{enumerate}
\item Up to isomorphism, the functors are invariants of oriented tangles, i.e. if
$t\cong t'$ then $F(t)\cong F(t')$.
\item In the Grothendieck group of the homotopy category of complexes of projective functors we have
the equality
\begin{eqnarray*}
q^k
\left[F_k\left(\;\parbox[thb]{0.4cm}{\includegraphics{cross1.eps}}\;\right)\right]
-q\left[{F}_k\left(\parbox[thb]{0.4cm}{\includegraphics{identity.eps}}\right)\right]
=q^{-k}\left[{F}_k\left(\;\parbox[thb]{0.4cm}{\includegraphics{cross2.eps}}\;\right)\right]
-q^{-1}\left[{F}_k\left(\parbox[thb]{0.4cm}{\includegraphics{identity.eps}}\right)\right],
\end{eqnarray*}
where $q^j$ means that the grading is shifted up by $j$.
\end{enumerate}
\end{theorem}

In other words, we get a categorification of the polynomial
$\mathfrak{sl}(k)$-invariant $\mathbf{P}_k$. Note that this is only a
categorification in the weak sense, which means we do not specify isomorphisms
defining the relations. This is somehow the drawback of our "down-to-earth"
combinatorial approach: we cannot control these morphisms.

In the last section, however, we bring the natural transformation into the
picture. For that we stick to the case $k=3$ as in \cite{Ksl3} (but see the
general Conjecture~\ref{conjecture}). To each basic foam as depicted in
Figure~\ref{basicfoams}, we associate just the obvious natural transformation
of functors given by adjointness properties. Now, any such natural
transformation defines a homomorphism when evaluating at any single object, in
particular if we evaluate it at the antidominant projective module in the most
regular block to choose from. Under Soergel's combinatorial functor $\mV$ this
morphism turns into a morphism between certain modules over the endomorphism
ring of the antidominant projective modules. These endomorphism rings have
however a very easy description, namely each of them is isomorphic to the
cohomology ring of some partial flag variety which are in most cases just
Grassmannians. Hence we finally end up with maps between modules over certain
cohomology rings, in fact with tensor products of certain cohomology rings.
These turn out to be exactly the maps in \cite{Ksl3}. In general these maps
should give rise exactly to the maps from \cite{MackaayStosic}. Putting dots on
a foam means in our approach nothing else than multiplication with an element
of the centre of (a direct summand) of the category categorifying the boundary
web.

In light of \cite{Brundan} and \cite{StSpringer} one might expect that not only
the partial flag varieties, but also Springer fibres and Spaltenstein
varieties, and the combinatorics of their cohomology rings should play a
crucial role in the complete picture.

\subsubsection*{Notation:} In the following we will abbreviate
$\otimes_\mC$ as $\otimes$.

\subsubsection*{Acknowledgments:} We would like to thank Christina
Cobbold and Wolfgang Soergel for useful discussions and comments.

\section{Trivalent coloured graphs and intertwiners}
\label{Section2}
 Throughout the whole paper we fix an integer $k\geq 2$ and
denote by $V$ the natural $k$-dimensional representation of the quantum group
$U_q(\mathfrak{sl}_k)$ with generic parameter $q$, and fix the standard basis
$v_i$, $1\leq i\leq k$, of $V$ (see \cite[5A.1]{Ja3}).
For $1\leq i\leq k$ we have the
fundamental weights $\omega_i$ with the corresponding irreducible
$\cU_q(\mathfrak{sl}_k)$-modules $\bigwedge^iV$.

For any $i,j\in\{1,2,\ldots k\}$ we have the exterior powers
$\bigwedge^i V$, $\bigwedge^j V$, $\bigwedge^{i+j} V$ together with
the intertwiner maps
\begin{eqnarray*}
\begin{array}[tbh]{lllllllll}
\displaystyle \pi_{i,j}^{i+j}:&\bigwedge^i V\otimes \bigwedge^j
V&\rightarrow&\bigwedge^{i+j} V& &\pi_{i+j}^{i,j}:
&\bigwedge^{i+j}V&\rightarrow&\bigwedge^i V\otimes \bigwedge^j V.
\end{array}
\end{eqnarray*}
For explicit formulae describing the intertwiners relevant in our context, we
refer to the next paragraph.

\renewcommand{\thefigure}{\arabic{figure}}
\setcounter{figure}{-1}
\begin{figure}[htb]
\begin{center}
 \includegraphics{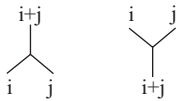}
\end{center}
  \caption{The graphs associated with $\pi_{i,j}^{i+j}$ and $\pi_{i+j}^{i,j}$ respectively}
  \label{fig:basicgraphs}
\end{figure}

There is a graphical description of intertwiners between tensor products of
exterior products of $V$ which associates to $\pi_{i,j}^{i+j}$ and
$\pi_{i+j}^{i,j}$ the coloured trivalent graphs as depicted in
Figure~\ref{fig:basicgraphs}. (Here and in the following the graphs should be
read from the bottom to the top.) Any arbitrary intertwiner can be described
via a composition of the elementary graphs from Figure~\ref{fig:basicgraphs},
so that one can associate with any intertwiner a trivalent graph coloured by
elements from the set $\{1,2,\ldots, k\}$ (which should be identified with the
set of fundamental weights for $\mathfrak{sl}_k$).


\subsection{Special intertwiners}
\label{specialint} In the context of knot and link invariants, a special role
is played by the pairs $(i,j)\in\{(1,1),(1,k-1),(k-1,1)\}$. We will use a (red)
very thick line for the labelling $k$, a  (green) thick line for the labelling
$k-1$. A (blue) normal line indicates the labelling by $2$, and finally a
thin black line indicates labelling by $1$. In the standard bases we have the
following explicit formulas:

\begin{eqnarray*}
\begin{array}[bt]{cccccccc}
{\includegraphics{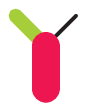}}:&\bigwedge^k
V&\rightarrow&\bigwedge^{k-1} V\otimes V\\
&w&\mapsto&\sum_{j=1}^k q^{j-1} w(j)\otimes v_j.\\
\includegraphics{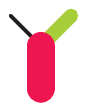}:&\bigwedge^k V&\rightarrow &V\otimes\bigwedge^{k-1} V\\
&w&\mapsto&\sum_{j=1}^k q^{k-j} v_j\otimes w(j).\\
\includegraphics{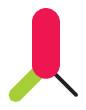}:&\bigwedge^{k-1} V\otimes V&\rightarrow&\bigwedge^{k}
V\\
&w(j)\otimes v_s&\mapsto&
\begin{cases}
q^{j-k} w& \text{ if $j=s$}\\
0&\text{ if $j\not=s$}
\end{cases}
\\
\includegraphics{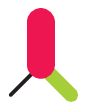}:&V\otimes\bigwedge^{k-1} V&\rightarrow&\bigwedge^{k}V\\
&v_s\otimes w(j)&\mapsto&
\begin{cases}
q^{j-1} w& \text{ if $j=s$}\\
0&\text{ if $j\not=s$}
\end{cases}\\
\includegraphics{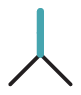}:&V\otimes V&\rightarrow&\bigwedge^2V\\
&v_i\otimes v_j&\mapsto&
\begin{cases}
q^{-1} v_i\wedge v_j& \text{ if $i>j$}\\
v_i\wedge v_j&\text{ if $i<j$}
\end{cases}\\
\includegraphics{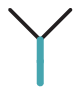}:&\bigwedge^2 V&\rightarrow &V\otimes V\\
&v_i\wedge v_j&\mapsto&v_j\otimes v_i+q v_i\otimes v_j\text{ if
$i<j$}
\end{array}
\end{eqnarray*}
where $w:=v_1\wedge v_2\wedge\ldots\wedge v_k$ and $w(j):=v_1\wedge\ldots
v_{j-1}\wedge v_{j+1}\wedge\ldots\wedge v_k$.

The relations between the intertwiners translate into relations between
trivalent graphs. Some of them - namely the ones involving only the special
intertwiners with labels from $\{1,2,k-1,k\}$ are depicted in the Relations (I)
to (IV) below.

These are the relevant graphs used in \cite{MOY} to define the
$\mathfrak{sl}_k$-invariants of links. Theorem~\ref{main1} gives a categorical
interpretation of these relations, including a functor valued
$\mathfrak{sl}_k$-invariant which enriches the polynomial invariant
$\mathbf{P}_k$.

\begin{figure}[htb]
  \centering
 \includegraphics{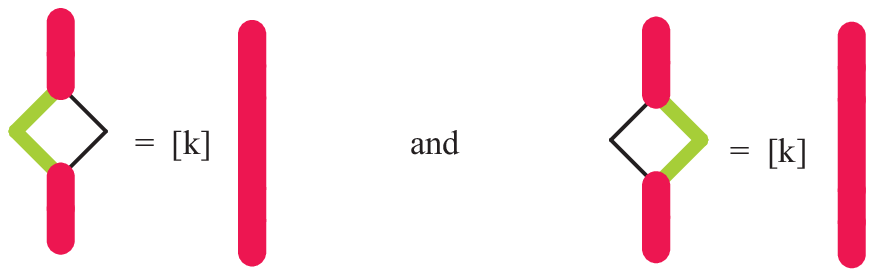}
  \caption{{\bf Relations (I)}: Two pairs of Intertwiners $\wedge^k V\rightarrow\wedge^kV$.}
  \label{fig:Rel1}
\end{figure}

\begin{figure}[htb]
  \centering
 \includegraphics{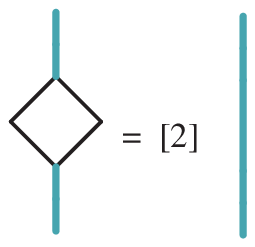}
  \caption{{\bf Relation (II)}: Intertwiners $\wedge^2 V\rightarrow\wedge^2V$}
  \label{fig:Rel2}
\end{figure}

\begin{figure}[htb]
  \centering
 \includegraphics{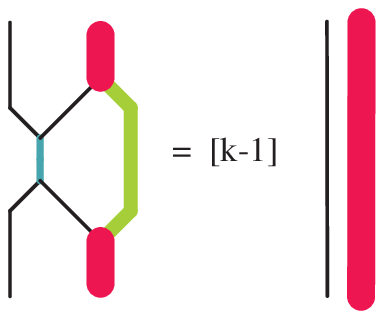}
  \caption{{\bf Relation (III)}: Intertwiners $V\otimes \wedge^k V\rightarrow V\otimes\wedge^kV$.}
  \label{fig:Rel3}
\end{figure}

\begin{figure}
  \centering
 \includegraphics{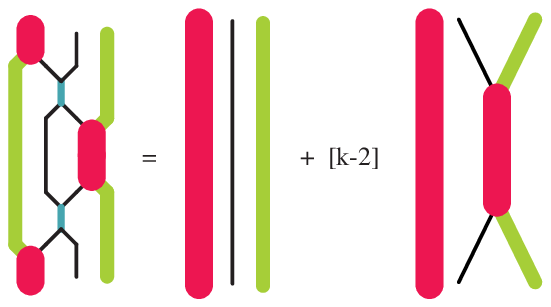}
  \caption{{\bf Relation (IV)}: Intertwiners $\wedge^k V\otimes V\otimes\wedge^{k-1} V\rightarrow\wedge^k V\otimes V\otimes\wedge^{k-1}V$.}
  \label{fig:Rel4}
\end{figure}

\begin{figure}
  \centering
 \includegraphics{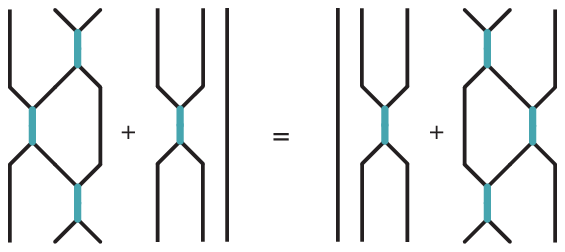}
  \caption{{\bf Relation (V)}: Intertwiners $V\otimes V\otimes V\rightarrow V\otimes V\otimes V$.}
  \label{fig:Rel5}
\end{figure}

\section{Box diagrams and fillings}
\label{Section3}
Fix a positive integers $n$. Any tensor product $V^{\otimes
i}$, exterior product $\wedge^i V$,  or combination of both, comes along with
the {\it standard basis} given by tensors of basis vectors of $V$ and exterior
products $v_{i_1}\wedge v_{i_2}\wedge\ldots\wedge v_{i_k}$ with strictly
decreasing indices $i_1>i_2>\ldots >i_k$.

A tuple $\mu=(\mu_1,\mu_2,\ldots,\mu_l)$ of nonnegative integers with
$\sum_{i=1}^l\mu_i=n$ is a {\it composition} of $n$, denoted $\mu\vDash n$. We
call the number $l$ the {\it length} ${l}(\mu)$ of $\mu$, and the number
 of non-zero entries of $\mu$ the {\it actual length}, denoted  ${ll}(\mu)$, of $\mu$. Associated with any composition $\mu$
 we have the {\it box diagram} $D^\mu$  - drawn in the $\{(x,y)\mid x\geq0, y\leq0\}$-quadrant of
the plane, numbering rows $1, 2,\ldots$ from top to bottom and columns $1,
2,\ldots$ from left to right - with $\mu_i$ boxes in the $i$th-column, placed
in the rows $1$ to $\mu_i$ - see Examples below.

Given a box diagram $D^\mu$ of type $\mu$ and a second composition
$\nu$ of $n$, a {\it filling} of $D^\mu$ {\it of type $\nu$} is a
filling of $D^\mu$ such that for $1\leq i\leq l(\nu)$, the number
$i$ appears exactly $\nu_i$ times. The filling is {\it column
strict} if in each column the numbers are strictly increasing from
top to bottom. If $l(\mu)\leq k$ we associate to a given column
strict filling $F$ of type $\nu$ of $D^\mu$ a standard basis element
$$\Phi(F)\in\bigwedge^\nu
V:=\bigwedge^{\nu_1}{V}\otimes\bigwedge^{\nu_2}{V}\otimes\ldots\otimes\bigwedge^{\nu_k}{V}$$
as follows: Let $c_{i,1}<c_{i,2}<\ldots <c_{i,\nu_i}$ be the numbers of the
columns, where the entry $i$ occurs, then
\begin{eqnarray}
\label{Phi} \Phi(F):=w_1\otimes w_2\otimes \ldots\otimes w_k
\end{eqnarray}
where $w_{i}:=v_{c_{i,1}}\wedge v_{c_{i,2}}\wedge \ldots\wedge
v_{c_{i,\nu_i}}$.

\begin{exs}
\label{Exbox} {\rm Let $n=6$, $k=3$, $\nu=(2,3,1)$. Then $\bigwedge^\nu V$ has
dimension $9$. For $\mu$ equal to $(3,2,1)$, $(3,1,2)$, $(2,1,3)$, $(2,3,1)$,
$(1,2,3)$, $(1,3,2)$ there is only one possible column strict filling of type
$\nu$ giving rise to the following basis vectors
\begin{eqnarray*}
(3,2,1)&\rightsquigarrow&v_1\wedge v_2\otimes v_1\wedge v_2\wedge v_3\otimes v_1\\
(3,1,2)&\rightsquigarrow&v_1\wedge v_3\otimes v_1\wedge v_2\wedge v_3\otimes v_3\\
(2,1,3)&\rightsquigarrow&v_2\wedge v_3\otimes v_1\wedge v_2\wedge
v_3\otimes v_3\\
(2,3,1)&\rightsquigarrow&v_1\wedge v_3\otimes v_1\wedge
v_2\wedge v_3\otimes v_1,\\
(1,2,3)&\rightsquigarrow&v_1\wedge v_2\otimes v_1\wedge v_2\wedge v_3\otimes v_2,\\
(1,3,2)&\rightsquigarrow&v_2\wedge v_3\otimes v_1\wedge v_2\wedge
v_3\otimes v_2.
\end{eqnarray*}

 For
$\mu=(2,2,2)$ there are the following three possible column strict
fillings with corresponding basis vectors
\small
\begin{eqnarray*}
\begin{array}{cccc}
\young(2,3)\young(1,2)\young(1,2)&
\young(1,2)\young(2,3)\young(1,2)&
\young(1,2)\young(1,2)\young(2,3)\\
v_2\wedge v_3\otimes v_1\wedge v_2\wedge v_3\otimes v_1& v_1\wedge
v_3\otimes v_1\wedge v_2\wedge v_3\otimes v_2& v_1\wedge v_2\otimes
v_1\wedge v_2\wedge v_3\otimes v_3.
\end{array}
\end{eqnarray*}
\normalsize
Let $n=2$, $k=3$, $\nu=(1,1)$, hence $\bigwedge^\nu V=V\otimes V$. Then we have
for instance the following box diagrams, where the dots are indicating the
columns with no boxes:
$$D^{(2,0,0)}=
\parbox[th]{0.5cm}
{\yng(1,1)}\bullet\;\bullet\quad\quad\quad D^{(0,2,0)}=
\bullet\;\;\parbox[th]{0.5cm}{\yng(1,1)}\bullet\quad\quad\quad
D^{(0,2,0)}=\bullet\;\bullet\yng(1,1)
$$
In each case there is only one possible column strict filling of type
$\nu=(1,1)$, namely {\parbox[h]{0.5cm}{\young(1,2)}}. The corresponding basis
elements of $V\otimes V$ are then $v_1\otimes v_1$, $v_2\otimes v_2$ and
$v_3\otimes v_3$ respectively. Figuring out the remaining basis vectors is left
to the reader. }
\end{exs}

\subsection{Actions of the symmetric group}
Let $\hat{\mathbf D}^\mu_\nu$ (resp. ${\mathbf D}^\mu_\nu$) be the
set of box diagrams of type $\mu$ with fillings (resp. column strict
filling) of type $\nu$. If $\nu=(1^n):=(1,1,\ldots 1)$ we will normally
omit the index $\nu$ in the notation. There is a special element
$T^\mu\in{\mathbf D}^\mu$ with the {\it standard filling} given by
putting the numbers $1,2,3,\dots, n$ in this order column by column
from the top to the bottom; for instance
$$T^{(2,2,2)}=\parbox[th]{5cm}{\young(1,2)\young(3,4)\young(5,6)}.$$
The {\it $i$-th box} of $D^\mu$ is the box with the number $i$ in
the standard filling; it is denoted by $b_i(D^\mu)$. Let $S_n$ be
the symmetric group with the usual generators $s_i$, $1\leq i\leq
n-1$. Then $S_n$ acts on $\hat{\mathbf D}^\mu_\nu$ from the right by
permuting the entries and from the left by permuting the boxes (with
their entries).
\begin{ex}
{\rm $T^{(2,2,2)}s_4s_3s_2=\young(1,3)\young(4,5)\young(2,6)$,
whereas $s_2s_3s_4T^{(2,2,2)}=\young(1,5)\young(2,3)\young(4,6)$.}
\end{ex}

\subsection{The correspondence}\label{sec3.2}
For any composition $\mu$ of $n$ let $\tilde{\mu}$ be the {\it reduced
composition} obtained by disregarding the zero entries of $\mu$. Let $S_\mu$ be
the corresponding Young subgroup, i.e. $S_\mu=S_{\tilde\mu_1}\times
S_{\tilde\mu_2}\times S_{\tilde\mu_{ll(\mu)}}$ of $S_n$. We denote by ${}^\mu
S_n$ the set of shortest coset representatives in $S_\mu\backslash S_n$,
similarly let $S_n^\mu $ be the set of shortest coset representatives in
$S_n/S_\mu$. Let $\bf{O}^\mu_\nu$ denote the set of cosets $c\in S_n/S_\nu$
such that $w\in {}^\mu S_n$ for any $w\in c$.

Assume we have a box diagram $D$ and $\nu\vDash n$. Then any filling of type
$\nu$ can be transferred into a filling of type $(1^n)$ by replacing first all
ones by the numbers $1,2,\ldots, \nu_i$ from left to right, then all two's by
the numbers $\nu_1+1, \ldots \nu_1+\nu_2$ etc. On the other hand, if we have a
filling $F$ of type $(1^n)$ then we can replace the first $\nu_1$ numbers by
$1$'s, the next $\nu_2$ numbers by $2$'s etc. We call the result $\psi_\nu(F)$.
The latter is always an element of $\hat{{\mathbf D}}^\mu_\nu$, but not
necessarily of ${\mathbf D}^\mu_\nu$. We have however the following result

\begin{lemma}
\label{bijections}
\begin{enumerate}
\item Let $\mu\vDash n$ and $w\in S_n$, then $wT^\mu=T^\mu w$.
\item The map $\Phi$ from \eqref{Phi} defines a bijection
\begin{eqnarray*}
\Phi_\nu^\mu:\quad\bigcup_{l(\mu)\leq k}{\mathbf D}^\mu_\nu&\oneone&\text{
elements of the standard basis of  } \bigwedge^\nu {V}
\end{eqnarray*}
\item There is a bijection
\begin{eqnarray*}
\Psi_\nu^\mu: \quad\bf{O}^\mu_\nu&\oneone&{\mathbf D}^\mu_\nu\\
w&\mapsto&\psi_\nu(wT^\mu).
\end{eqnarray*}
\end{enumerate}
\end{lemma}

\begin{proof}
By definition, the entry of the $i$th-box of $T^\mu$ is precisely $i$, so the
first statement is obvious. The map $\Phi_\nu^\mu$ is obviously injective. To
see that it is surjective note that a basis of $\bigwedge^\nu V$ is given by
elements of the form $w_1\otimes w_2\otimes \ldots\otimes w_k$ where
$w_{i}:=v_{c_{i,1}}\wedge v_{c_{i,2}}\wedge \ldots\wedge v_{c_{i,r_i}}$, where
for any $i$ we have $c_{i,1}<{c_{i,2}}<\ldots <{c_{i,r_i}}$ and $1\leq
c_{i,j}\leq k$. A preimage of $w_1\otimes w_2\otimes \ldots\otimes w_k$ can be
constructed as follows: we create a box diagram with column strict filling by
putting ones in the columns $c_{(1,j)}$, then $2$'s in the columns $c_{(2,j)}$
etc. As a result we get an element in $\bigcup_{l(\mu)\leq k}{\mathbf
D}^\mu_\nu$ which is obviously a preimage, and $\Phi_\nu^\mu$ is surjective.

Let's take the box diagram $T^\mu$ associated with $\mu$ with the standard
filling. $S_n$ acts transitively from the left on $\hat{\mathbf D}^\mu$ giving
rise to a bijection $\alpha: S_n\cong S_n T^\mu$. From the definition of the
left action of $S_n$ on diagrams with fillings we get directly that $wT^\mu \in
{}^\mu S_n$ if and only if, in each column, the entries are strictly increasing
from top to bottom. Hence $\Psi_\nu^\mu$ is a bijection if $\nu=(1^n)$. If
$\nu$ is now arbitrary, then $w\in\bf{O}^\mu_\nu$ if and only if the entries in
the columns are still strictly increasing from top to bottom if we replace the
first $\nu_1$ numbers by ones, then the next $\nu_2$ numbers by twos etc. The
claim is then obvious.
\end{proof}

For any set $M$ let $\mC[M]$ be the free $\mC$-module with basis given by the
elements of $M$. If $\mu=1^n$, then the action of $S_n$ turns $\mC[{\mathbf
D}^\mu]$ into the permutation module, which is a special instance of the
induced sign module $N(\mu)=\mC[S_n]\otimes_{\mC[S_\mu]}\operatorname{sgn}$ for
arbitrary $\mu$. The latter has a basis given by $x\otimes 1$, $x\in {}^\mu
S_n$. We identify this space in the obvious way with $\mC[\bf{O}^\mu]$ and
$\mC[{\mathbf D}^\mu]$ so that Lemma~\ref{bijections} induces isomorphisms of
$S_n$-modules $\mC[\bf{O}^\mu]\cong\mC[{\mathbf D}^\mu]$ and
$\bigoplus_\mu\mC[{\mathbf D}^\mu]\cong {V}^{\otimes n}$ where $S_n$ acts by
permuting the factors.

All this can be quantised: If ${}^\mZ\bf{O}^\mu_\nu$ denotes the free
$\mC[q,q^{-1}]$-module with basis $\bf{O}^\mu_\nu$ then we view
${}^\mZ\bf{O}^\mu$ as the induced sign module
$\cN(\mu)=H_q(S_n)\otimes_{H_q(S_\mu)}\operatorname{sgn}$ for the Iwahori-Hecke
algebra $H_q(S_n)$ and
$$\bigoplus_\mu{}^\mZ{\bf O}^\mu\cong V^{\otimes n}$$ where
$H_q(S_n)$ acts via the $R$-matrix.

The Hecke algebra $H_q(S_n)$ comes along with the standard basis $H_x$, $x\in
S_n$, and with the Kazhdan-Lusztig basis $\underline H_x$, $x\in S_n$. In the
following we use the normalisation of \cite{SoKipp}. In particular, $\underline
H_s=H_s+qH_e=:H_s+q$. Associated with $x\in S_n$ we have $(t(x), t'(x))$, the
corresponding pair of standard tableaux via the Robinson-Schensted
correspondence. We will need the following well-known result (see e.g.
\cite[Section3]{Martin}): If $t(x)$ has more than $ll(\mu)$ rows then
$\underline{H}_x$ is in the annihilator of ${}^\mZ\bf{O}^\mu$.

\subsection{Category $\cO$} \label{categoryO}
We consider the Lie algebra $\mathfrak{gl}_n$ and the
corresponding Bern\-stein-Gelfand-Gelfand category
$\cO=\cO(n)$ associated with the
standard triangular decomposition $\mathfrak{gl}_n=\mathfrak{n}_-\oplus\mathfrak{\mh}\oplus
\mathfrak{\mn}=\mathfrak{n}_-\oplus\mathfrak{b}$,
see \cite{BGG}. The Weyl group is identified with the permutation group $S_n$
in the standard way.

For any composition $\la$ of $n$ we {\bf fix} an integral block
$\cO_{\tilde\la}$ of $\cO$ such that the projective Verma module in this block
has highest weight $\tilde{\la}$, and the stabiliser of $\tilde{\la}$ is
$S_\la$. By abuse of notation we denote this block by $\cO_\la$ and the highest
weight of the projective Verma module
$P(\tilde\la)=M(\tilde\la)\in\cO_{\tilde\la}$ by $\la$. For $\mu\vDash n$ let
$\cO^\mu_\la$ be the subcategory given by all locally $\p$-finite objects,
where $\p$ is the parabolic (containing $\mb$) with Weyl group $S_\mu$. The
simple objects in $\cO_\la$ are exactly the simple highest weight modules
$L(x\cdot\la)$ with $x\in S_n^\la$, with the corresponding Verma modules
$M(x\cdot\la)$. The simple objects in $\cO_\la^\mu$ are exactly the simple
highest weight modules $L(x\cdot\la)$ with $x\in \bf{O}_\la^\mu$. In
particular, $\mC[\bf{O}_\la^\mu]$ can be identified with the complexified
Grothendieck group of $\cO_\la^\mu$ by mapping $x\in\bf{O}_\la^\mu$ to the
isomorphism class of the parabolic Verma module $M^\mu(x\cdot\la)$ with highest
weight $x\cdot\la$.

We denote by $^\mZ\cO_\la^\mu$ the graded version of $\cO_\la^\mu$ as
introduced in \cite{BGS} and further developed in \cite{Stgrad} and
\cite[Section 2]{StDuke}. Each parabolic Verma module $M^\mu(x\cdot\la)\in
\cO_\la^\mu$ has a standard graded lift $\Delta(x\cdot\la)\in{}^\mZ\cO_\la^\mu$
with head concentrated in degree zero. For $j\in\mZ$ we denote by
$\Delta(x\cdot\la)\langle j\rangle$ the lift with head in degree $j$, in
particular $\Delta(x\cdot\la)\langle 0\rangle=\Delta(x\cdot\la)$. Let
$P(x\cdot\la)\langle j\rangle$ be the indecomposable projective cover
of $\Delta(x\cdot\la)\langle j\rangle$. More
generally, we denote by $\langle j\rangle$ the functor which shifts the grading
up by $j\in\mZ$.

Note that the complexified Grothendieck group $[^\mZ\cO_\la^\mu]$ of
$^\mZ\cO_\la^\mu$ is naturally isomorphic to $^\mZ\bf{O}_\la^\mu$ by mapping
$\Delta(x\cdot\la)\langle j\rangle$ to $q^j x$. In the following we will abuse
notation and denote $\Delta(x\cdot\la)\langle j\rangle$ by
$q^j\Delta(x\cdot\la)$ or even by $q^j\Delta(x)$ or $q^j \Delta(i_1\;i_2\ldots
i_r)$ if $x=s_{i_1}\ldots s_{i_r}$ is a reduced expression for $x$ and it is
clear from the context to which category the module belongs to. Analogous
abbreviations will be used for the projectives $P(x\cdot\la)\langle j\rangle$.

\section{The same combinatorics in three disguises}
\label{Section4}
\subsection{Translation functors - combinatorially}
\label{translations} We first recall the explicit combinatorics of special
projective functors, namely the translation functors on and out of the walls.
Thanks to {\bf Fact 1} below the combinatorics describes the functor
completely.

Let $\la,\mu\vDash n$. If $S_\la\subseteq S_\nu$ then there is the {\it
translation out of the walls} functor (see \cite[4.11]{Ja2})
$$T_{\nu}^{\lambda}:\;\cO(n)_{\nu}\longrightarrow
\cO(n)_{\lambda}$$ with its {\it standard graded lift}
$$\theta_{\nu}^{\lambda}:\;{}^\mZ\cO(n)_{\nu}\longrightarrow
{}^\mZ\cO(n)_{\lambda}$$ which is uniquely determined by requiring that
$\Delta(e)$ is mapped to a standard lift of the (indecomposable) projective
module $T_\nu^\la M(\nu)$. In the following we will only need special instances
of translation functors (analogous to our special choices of intertwiners in
Section~\ref{specialint}). Let $\nu, \lambda\vDash n$ such that there exists
some $l$ such that $\lambda_t=\nu_t$ for $t<l$, $\lambda_{t+1}=\nu_{t}$ for
$t>l+1$ and set $j=\lambda_1+...\lambda_{l-1}$.
\begin{enumerate}[({\text Case} 1.)]
\label{Condition}
\item If moreover $\lambda_l=1$, $\lambda_{l+1}=i$, $\nu_l=i+1$ then
$\theta_{\nu}^{\lambda}:\;{}^\mZ\cO(n)_{\nu}\longrightarrow
{}^\mZ\cO(n)_{\lambda}$ maps $\Delta(e)$ to the graded projective module
$P((i+j)(i+j-1)\ldots
(j+1))$. The latter has each of the following:
$$\Delta((j+i)...(j+1)),\; q\Delta((j+i-1)...(j+1)),\; ..., \;
q^i\Delta(e)$$ exactly
once as graded Verma subquotients. To abbreviate this we will say $\Delta(e)$
is mapped to $A_{j+i}^{j+1}$ as defined in \eqref{AB}.
\item If moreover $\lambda_l=i$, $\lambda_{l+1}=1$, $\nu_l=i+1$,
then $\theta_{\nu}^{\lambda}:\;{}^\mZ\cO(n)_{\nu}\longrightarrow
{}^\mZ\cO(n)_{\lambda}$ maps $\Delta(e)$ to the graded projective
$P\left((j+1)\;(j+2)\ldots (j+i)\right)$ which has
\begin{gather*}
 \Delta\left((1+j)\ldots (i+j)\right),\; q\Delta\left((2+j)\ldots (i+j)\right),\; \ldots,\\ q^{i-1}\Delta(i+j), \; q^{i}\Delta(e),
\end{gather*}
as graded Verma subquotients. In a short form we say that $\Delta(e)$ is
mapped to $B_{j+1}^{i+j}$ as defined in
\eqref{AB}.
\end{enumerate}

  Translation
functors preserve parabolic subcategories, hence it makes sense to
define
\begin{eqnarray*}
\bigcurlyvee_{(i+j)}^{(i,\,\,j)}&=&\bigoplus_\mu
\theta_{(i+j)}^{(i,j)}:\quad\bigoplus_\mu{}^\mZ\cO(i+j)^\mu_{(i+j)}\longrightarrow
\bigoplus_\mu{}^\mZ\cO(i+j)^\mu_{(i,j)}.
\end{eqnarray*}
where the sum runs over all compositions $\mu$ of length at most
$k$.

Again if we have $S_\la\subseteq S_\nu$ there is also the {\it translation onto
the walls} functor $T_{\la}^\nu:\cO(n)_\la\longrightarrow \cO(n)_\nu$.
We have its {\it standard graded lift}
\begin{eqnarray*}
\theta_{\la}^\nu:{}^\mZ\cO(n)_\la\longrightarrow {}^\mZ\cO(n)_\nu
\end{eqnarray*}
which maps $\Delta(x)$ to $q^{-r}\Delta(z)$, where $z$ and $r$ are defined by
writing $x=zy$ with $y\in S_{\nu}$, and $z\in S_n^\nu$ a shortest coset
representative and $r=l(y)$ being the length of $y$. We define
\begin{eqnarray*}
\bigcurlywedge^{(i+j)}_{(i,\,\,j)}=\bigoplus_\mu
\theta_{(i,j)}^{(i+j)}:\quad\bigoplus_\mu{}^\mZ\cO(i+j)^\mu_{(i,j)}
\longrightarrow
\bigoplus_\mu{}^\mZ\cO(i+j)^\mu_{(i+j)},
\end{eqnarray*}
where the sum runs over all compositions $\mu$ of length at most
$k$.

Let $\la$, $\nu$, $\mu$ be compositions of $n$. Translation functors out and
onto walls are special instances of projective functors. We denote by
$\cP(\la,\nu)$ the set of projective functors from $\cO_\la$ to $\cO_\nu$ as
introduced and classified in \cite{BG}. We recall the following well-known
facts:

\begin{enumerate}[{\bf{\text Fact} 1}]
\item (\cite{BG}) A projective functor $F\in\cP(\la,\nu)$ is (up to
isomorphism) completely determined by its value on $ M(\la)$, i.e. we have an
isomorphism of projective functors $F\cong G$ if and only if there is an
isomorphism of modules $F M(\la)\cong G M(\la)$. More precisely: $F
M(\la)\in\cO_\nu$ is projective and $F$ decomposes into indecomposable summands
exactly according to the decomposition of $F M(\la)$ into indecomposable direct
summands.

\item (\cite[Corollary 3.12]{StDuke}, \cite{Stgrad}) Let $F\in\cP(\la,\nu)$ be indecomposable. There exists a
graded lift $\tilde{F}:{}^\mZ\cO_\la\rightarrow{}^\mZ\cO_\nu$. Up to
isomorphism and shift in the grading it is unique, and up to isomorphism
completely determined by its value on $\Delta(e)\in {}^\mZ\cO_\la$ (thanks to
{\bf Fact 1}).

\item (\cite[Proposition 4.2]{StDuke} and references therein) Let $F\in\cP(\la,\nu)$ be indecomposable such that
$\theta_\nu^{(1^n)} F M(\la)\cong P(x)$. Assume the tableau $t(x)$ has more
than $k$ rows. Then the restriction of $F$ to $\cO_\la^\mu$ is zero for any
$\mu$ with $ll(\mu)\leq k$.
\end{enumerate}

\subsection{The combinatorial action of trivalent graphs}
We define $\mC[q,q^{-1}]$-linear maps

\begin{eqnarray*}
\bigcurlywedge_{(i,\;\;j)}^{(i+j)}:\quad\bigoplus_\mu\mC[{\mathbf D}_{(i,j)}^\mu]&\rightarrow&\bigoplus_\mu\mC[{\mathbf D}_{(i+j)}^\mu]\\
\bigcurlyvee^{(i,\;\;j)}_{(i+j)}:\quad\bigoplus_\mu\mC[{\mathbf
D}_{(i+j)}^\mu]&\rightarrow&\bigoplus_\mu\mC[{\mathbf D}_{(i,j)}^\mu]
\end{eqnarray*}
where $\mu$ runs over all compositions of $i+j$ with at most $k$ parts, as
follows: In the first case we write any box diagram $D$ with filling of type
$(i,j)$ as $D=x\psi_{(i,j)}T^\mu$ with $x$ of smallest possible length. Then
$D$ is mapped to a box diagram $q^{-l(x)}D'$ where $D'\in{\mathbf
D}^\mu_{(i+j)}$ has the same shape as $D$, but for the filling we replace the
$2$'s by $1$'s. In the second case a box diagram $D$ of type $\mu$ and filling
$(i+j)$ is mapped to $\sum_{I}q^{l_I}D_I$, where $I$ runs through all possible
 subsets of cardinality $j$ of the set of boxes of $D$. The diagram
$D_I$ is obtained from $D$ by replacing all $1$'s in the boxes from
$I$ by $2$'s, and $l_I$ is equal to $ij$ minus the length of the
element $x$ of minimal length such that $D_I=x\psi_{(i,j)}T^\mu$.

\begin{ex}
Let $\nu=(3)$, $\nu'=(2,1)$ and $r=1$. Then
\small
\begin{eqnarray*}
\bigcurlyvee^{(2,1)}_{(3)}(\young(1)\young(1)\young(1))=
&\young(2)\young(1)\young(1)+q\young(1)\young(2)\young(1)
+q^2\young(1)\young(1)\young(2)&\\
\bigcurlywedge_{(2,1)}^{(3)}(\young(2)\young(1)\young(1))
=q^{-2}\young(1)\young(1)\young(1)&\displaystyle
\bigcurlywedge_{(2,1)}^{(3)}(\young(1)\young(2)\young(1))
=q^{-1}\young(1)\young(1)\young(1)&
\bigcurlywedge_{(2,1)}^{(3)}(\young(1)\young(1)\young(2))
=\young(1)\young(1)\young(1)
\end{eqnarray*}
\end{ex}
\normalsize

We have the obvious generalisation of this procedure if $\la$ and $\nu$ are of
the form as in (Case 1) or (Case 2), namely the role played by the entries $1$
and $2$ above is then the role of $j+1$ and $j+2$. This defines the maps
$\displaystyle \bigcurlyvee^\la_\nu:\bigoplus_\mu\mC[{\mathbf
D}_{\nu}^\mu]\rightarrow\bigoplus_\mu\mC[{\mathbf D}_{\la}^\mu]$ and
$\displaystyle \bigcurlywedge_\la^\nu:\quad\bigoplus_\mu\mC[{\mathbf
D}_{\la}^\mu]\rightarrow\bigoplus_\mu\mC[{\mathbf D}_{\nu}^\mu]$, where $\mu$
runs always through all compositions of $n$ with at most $k$ parts.

\begin{prop}
For simplicity let $\la$ and $\nu$ be as in Case 1 or Case 2. The
following diagram commutes:
\begin{eqnarray*}
\xymatrix{ \left[\bigoplus_\mu {}^\mZ\cO_\la^\mu\right]\ar[d]^{[F]}
\ar[r]^{\Xi_\la}
&\bigoplus_\mu\mC[{\bf D}_\la^\mu]\ar[r]^{\Phi_\la}\ar[d]^{G}
&\bigwedge^\la V\ar[d]^{H}
\\
\left[\bigoplus_\mu{}^\mZ\cO_\nu^\mu\right] \ar[r]^{\Xi_\nu}
&\bigoplus_\mu\mC[{\bf D}_\nu^\mu] \ar[r]^{\Phi_\nu} &\bigwedge^\nu V
}
\end{eqnarray*}
where $F=\theta_\la^\nu$ is the standard lift of the translation functor
to the wall, $\displaystyle G=\bigcurlywedge_\la^\nu$, $H=\pi_\la^\nu$ is
the corresponding intertwiner, and the $\Phi$'s and the $\Xi$'s
are the maps given via
all the identifications described in Subsections~\ref{sec3.2} and
\ref{categoryO}. The analogous statement holds if the roles of $\la$
and $\nu $ are swapped.
\end{prop}

\begin{proof}
The proof is a straightforward checking and therefore omitted.
\end{proof}

\section{Functor-valued invariants of coloured trivalent graphs}
\label{Section5} In this section we will indicate how to construct a
functor-valued invariant of trivalent graphs. Since we are mainly interested in
invariants of knots, we stick to what we called the special intertwiners
together with the Relations (I) to (V).

For a basic trivalent graph as depicted in Figure~\ref{fig:basicgraphs} we
associate the corresponding translation functor from
Section~\ref{translations}, more precisely let $\la\vDash n$ and $\nu\vDash m$
and assume we have a basic intertwiner $\bigwedge^\nu V\rightarrow\bigwedge^\la
V$ or its corresponding graph. Then we first associate as an intermediate step
the corresponding {\it non-parabolic translation functor}
$\theta_\nu^\la:{}^\mZ\cO(m)_\la\rightarrow{}^\mZ\cO(n)_\nu$ and call it the
{\it naively associated functor}. Afterwards we take the direct sum of all the
restriction to all parabolic with at most $k$ parts. The result is what we call
{\it the functor associated with the intertwiner} or {\it the functor
associated with the graph} we started with.

We will need the following

\begin{enumerate}[{\bf Fact 4}]
\item Let $F:\;_{}^\mZ\cO_\la\rightarrow\; _{}^\mZ\cO_\la$ be a
composition of functors naively associated to any of the graphs depicted in
Relation (I) to Relation (IV). Then we have
$F\Delta(\la)\cong P$ where $P$ is a finite
direct sum of graded projective modules from the set
$$\{Q, Q\langle k\rangle \oplus Q\langle -k\rangle\mid k\in\mZ\}$$
where $Q$ runs through the standard lifts of indecomposable projective module
in $\cO_\la$.
\end{enumerate}

\begin{proof}
Let $d$ be the usual duality on $\cO$. Let $F'=T_\la^\nu$ be a translation on
or out of the walls with $\la$ and $\nu$ related as in (Case 1) or (Case 2).
Then $dF'\cong F'd$ (\cite[4.12(9)]{Ja2}). Let $\operatorname{d}$ be the
standard graded lift of the duality (\cite[6.1.1]{Stgrad}). An easy direct
calculation shows that $\operatorname{d}F'\cong F'\operatorname{d}\langle
2(n_\nu-n_\la)\rangle$, where $n_\nu$ (resp. $n_\la$) denotes the length of the
longest element in $S_\nu$ (resp. $S_\la$). In particular,
$\operatorname{d}F\cong F\operatorname{d}$. Let $\operatorname{T}$ be the
graded lift of the twisting functor (\cite{AS}, \cite[Section 5]{FKS})
corresponding to the longest element $w_0$ of the Weyl group such that
$\operatorname{T}\Delta(x\cdot\la)$ is mapped to
$\operatorname{d}\Delta(w_0x\cdot\la)$. Let first $\theta=\theta_\la^\nu$ be a
translation onto the walls with $\la$ and $\nu$ related as in (Case 1) or (Case
2). Then $\operatorname{T}\theta=\theta\operatorname{T}$ if we forget the
grading (\cite{AS}), and then
$\operatorname{T}\theta=\theta\operatorname{T}\langle r\rangle$ for some
integer $r$.

Analogously, $\operatorname{T}\theta'=\theta'\operatorname{T}\langle s\rangle$
for some $s\in\mZ$, where $\theta'=\theta^\la_\nu$. Hence,
$\theta\theta'\operatorname T\langle s\rangle=\theta\operatorname
T\theta'=\operatorname{T}\theta\theta'\langle -r\rangle$. Since $\theta\theta'$
is just the direct sum of several copies of the identity functors (possibly
shifted in the grading), we get $s=-r$. Since all the functors to consider are
associated with graphs having a reflection symmetry in a vertical line, the sum
of overall shifts is zero. This means $\operatorname{d}F\Delta(\la)\cong
F\operatorname{d}\Delta(\la)\cong
F\operatorname{T}^2\Delta(\la)\cong\operatorname{T}^2 F\Delta(\la)$, and since
$\operatorname{d}$ maps $Q\langle k\rangle$ to $(\operatorname{d}Q)\langle
-k\rangle$, whereas $\operatorname{T}^2$ maps $Q\langle k\rangle$ to
$(\operatorname{d}Q)\langle k\rangle$ (see \cite[Section 3]{AS}) the statement
follows.
\end{proof}

Let us summarise what we have: we associated to each trivalent graph two
functors the naively associated one and then the direct sum of its restriction
to all parabolics attached to a composition with at most $k$ parts. We will
show that the latter functors satisfy the Relations (I) to (V). Thanks to {\bf
Fact 1} to {\bf Fact 4} this becomes a purely combinatorial problem, which also
shows that it is enough to verify the the relations of the functors locally,
without paying attention how complicated the graphs might be outside this small
region.

For any positive integers $r\geq s$, we will use the following abbreviations
\begin{eqnarray}
\label{AB}
\begin{array}[htb]{crrrlcrrr}
A_r^s&=&r\;(r-1)\;(r-2)\ldots s&&B_s^r&=&s\;(s+1)\;\ldots(r-1)\;r\\
&&q(r-1)\;(r-2)\ldots s&&&&q(s+1)\;\ldots (r-1)\;r\\
&&\vdots&&&&\vdots\\
&&q^{r-s}s&&&&q^{r-s}r\\
&&q^{r-s+1}e&&&&q^{r-s+1}e
\end{array}
\end{eqnarray}

In the following we will also ``multiply''
such (unordered) lists and write $A B$
to denote the list of all concatenations $ab$, where $a\in A$ and $b\in B$. For
instance, $A_2^1 B_3^4$ denotes the list $2\,1\,3\,4, q\,2\,1\,3, q^2\,2\,1,
q1\,3\,4, q^2\,1\,3, q^3 1, q^2\,3\,4, q^3\,3, q^4 e$.

We denote by $[m]=\frac{q^{m}-q^{-m}}{q-q^{-1}}$ the $m$-th quantum number.
For a list $A$ as above we denote by $[m]A$ the list containing
$q^{m-1-2j}a$, $0\leq j\leq m-1$ , $a\in A$.

For a basic trivalent graph as depicted in Figure~\ref{fig:basicgraphs} we
associate the corresponding translation functor from
Section~\ref{translations}. We are going to show now that the Relations (I) to
(IV) are satisfied. As a consequence we will obtain Theorem 1 from the
Introduction.

\begin{prop}[Relations (I) and (II)]
\label{propRel2}
Let $l\in\{k,2\}$. There are isomorphisms of
functors
$$\theta_{l-1,1}^{l}\;\theta^{l-1,1}_{l}\cong[l]\operatorname{id}\quad{\text
and}\quad \theta_{1,l-1}^{l}\;\theta^{1,l-1}_{l}\cong[l]\operatorname{id}.$$
Hence, the relations from Figures~\ref{fig:Rel1} and \ref{fig:Rel2} hold (even
for the naively associated functors).
\end{prop}

\begin{proof}

Thanks to {\bf Fact 2} it is enough to compare the image (even its Verma flag!)
of the functors applied to the projective Verma module $\Delta(e)$. The first
functor is going from the block with singularity $\nu=(l)$ to $\nu=(1,l-1)$ and
back to $(l)$. Combinatorially, the image of $\Delta(e)$ is given as follows:
\begin{eqnarray*}
\begin{array}{c|c|c}
(l)&(1,k-1)&(l)\\
\hline e&A_l^1&[l]e
\end{array}
\end{eqnarray*}
Here, the first row indicates the singularity $\nu$, whereas the second row
displays the Verma flag of the corresponding functor applied to $\Delta(e)$
according to the combinatorics of translation functors. The first isomorphism
follows then directly, the second is completely analogous. In particular, the
Relations (I) and (II) hold for both, the naively associated functors as well
as their parabolic versions.
(Note that our argument doesn't make any assumptions on $l$, hence
the statement is true in bigger generality.)
\end{proof}

\begin{prop}[Relation (III)]
\label{PropRel3} Let $G$ be the naively associated functor to the left hand
side diagram of Figure~\ref{fig:Rel3}. Then there is an isomorphism of functors
\begin{eqnarray}
\label{split}
G:=\theta_{(1,1,k-1)}^{(1,k)}\;\theta^{(1,1,k-1)}_{(2,k-1)}\;\theta_{(1,1,k-1)}^{(2,k-1)}\;\theta_{(1,k)}^{(1,1,k-1)}\cong
F\oplus [k-1]\operatorname{id},
\end{eqnarray}
where $F$ is indecomposable and vanishes when restricted to any parabolic with
at most $k$ parts. In particular, the relation depicted in Figure~\ref{fig:Rel3} holds.
\end{prop}

\begin{proof}
Combinatorially, the naively associated functor is given as follows:
\begin{eqnarray*}
\begin{array}{l|l|l|l|l}
   (1,k)&(1,1,k-1)&(2,k-1)&(1,1,k-1)     &(1,k)\\
\hline e&A_k^2   &A_k^2 &A_k^21, qA_k^2&A_k^1, [k-1]e
\end{array}
\end{eqnarray*}
Using {\bf Fact 4} we get that $G\cong G'\oplus [k-1]\operatorname{id}$, where
$G'$ maps $\Delta(e)$ to $P(k(k-1)\ldots 1)$. Now we use {\bf Fact 3} and
consider $\theta_{(1,k)}^{(1^{k+1})}G'\Delta(e)=P(x)$, where $x$ is the
following permutation (of $n=k+1$ letters)
\begin{eqnarray*}
x=\left(
\begin{array}[tbh]{cccccccccccccc}
1&2&3&\ldots&\ldots&k&k+1\\
k+1&k&k-1&k-2&\ldots&2&1
\end{array}
\right)
\end{eqnarray*}
Under the Robinson-Schensted algorithm this corresponds to a tableau with
entries $1,2,\ldots, k,k+1$ in its first column, hence has $k+1$ rows. By {\bf
Fact 3}, the functor $F$ is zero when restricted to any parabolic with at most
$k$ parts. Hence the statement follows.
\end{proof}

We also have to check the relation which we obtain by reflecting the graphs
from Figure~\ref{fig:Rel3} in a vertical line passing between the two graphs.
This can be done completely analogously as above. Alternatively, consider the
isomorphism of the Lie algebra $\mathfrak{gl}_n$ given by the obvious
involution of the Dynkin diagram which swaps the $i$-th with the $n-i$-th node.
This isomorphism defines an auto-equivalence of the category $\cO$ for
$\mathfrak{gl}_n$ which identifies $\cO(n)_\nu^\mu$ with
$\cO(n)_{\tilde\nu}^{\tilde{\mu}}$, where the partition are 'reflected in a
vertical line'. Applying this involution we are back at the situation described
in Proposition~\ref{PropRel3}.

\begin{prop}[Relation (IV)]
\label{propRel4} Let $G_3$ be the functor naively associated with the graph on
the LHS of Figure~\ref{fig:Rel4}. There is an isomorphism of functors
\begin{eqnarray*}
G_3\cong
F\oplus[k-2]\;\theta^{(k,1,k-1)}_{(k,k)}\theta_{(k,1,k-1)}^{(k,k)}\oplus
\operatorname{id}_{(k,1,k-1)}
\end{eqnarray*}
where $F$ is an indecomposable functor which vanishes when
restricted to any parabolic with at most $k$ parts. In particular,
the Relation displayed in Figure~\ref{fig:Rel4} holds.
\end{prop}

\begin{proof}
The functor $G_3$ is a composition of different translation
functors. We go, step by step, through the combinatorics:
\begin{eqnarray*}
\begin{array}{l|l|l|l}
   (k,1,k-1)&(k-1,1,1,k-1)&(k-1,2,k-1)&(k-1,1,1,k-1)\\
\hline e&B_1^{k-1}&B_1^{k-1}&B_1^{k-1} k,\; q B_1^{k-1}
\end{array}
\end{eqnarray*}
If we now go to $(k-1,1,k)$ nothing changes and back to
$(k-1,1,1,k-1)$ we obtain \small
\begin{eqnarray*}
\begin{array}{ll}
B_1^{k-1}k\;(2k-1)\;(2k-2)\ldots (k+1)\\
q B_1^{k-1}k\;(2k-2)\ldots (k+1)&q B_1^{k-1}\;(2k-1)\;(2k-2)\ldots (k+1)\\
q^2 B_1^{k-1}k\;(2k-3)\ldots (k+1)&q^2 B_1^{k-1}\;(2k-2)\ldots (k+1)\\
\vdots&\vdots\\
q^{k-1} B_1^{k-1}k&q^{k-1} B_1^{k-1}(k+1)\\
&q^{k} B_1^{k-1}
\end{array}\nonumber\\
\end{eqnarray*}
\normalsize We denote the column on the left hand side by $C_1$ and the one on
the right hand side by $C_2$ and define $D$ to be the $C_1$ where we remove the
part $q^{k-1}B_1^{k-1}k$. (i.e. all the graded Verma modules indexed by the
elements which become shorter if we multiply with $k$ from the right hand
side.) Note also that $C_1=B_1^{k-1}kA_{2k-1}^{k+1}$ and
$C_2=qB_1^{k-1}A_{2k-1}^{k+1}$.

If we pass from $(k-1,1,1,k-1)$ to $(k-1,2,k-1)$ and go back to
$(k-1,1,1,k-1)$ then our $C_1$ together with $C_2$ from above is
then turned into the collection
$$Dk,\, q D,\, C_2k,\, qC_2,\, q^{k-2}B_1^{k-1}k,\, q^{k-1}B_1^{k-1}.$$
 Finally, we have
to go to $(k,1,k-1)$. The elements $D k$, $qD$, $q^{k-2}B_1^{k-1}k$ and $C_2 k$
stay the same, $qC_2$ becomes $q\sum_{j=1}^{k}q^jq^{-(k-j)}A_{2k-1}^{k+1}$, and
$q^{k-1}B_1^{k-1}$ becomes $(1+q^2+q^4+\cdots q^{2(k-1)})e$. Together with {\bf
Fact 4}, we finally obtain the following decomposition into indecomposable
projective modules:
$$P(1\ldots k-1\;k\;(2k-1)\ldots (k+1))\oplus P(e)\oplus [k-2]P((2k-1)\;(2k-2)\ldots k+1).$$
Now it's time again to use {\bf Fact 3}: take the element $y=1\ldots
(k-1)k(2k-1)\ldots k$ and translate $P(y)$ out of all walls. We get
$P(yz)$, where $z$ is the longest element of $S_k\times S_1\times S_{k-1}$. Now
we write $yz$ as a permutation $x$ (of $n=2k$ letters),
\begin{eqnarray*}
x=\left(
\begin{array}[tbh]{cccccccccccccc}
1&2&3&\ldots&k&k+1&k+2&\ldots&n-1&n\\
k+1&k&k-1&\ldots&2&n&n-1&\ldots&k+2&1
\end{array}
\right)
\end{eqnarray*}
Under the Robinson-Schensted algorithm, $x$ corresponds to a tableau with
entries $1,2,\ldots, k,k+1$ in its first column, hence has $k+1$ rows.
Therefore, the functor $F$ is zero when restricted to any parabolic with at
most $k$ parts.
\end{proof}

The Relation from Figure~\ref{fig:Rel5} is nothing else than the Hecke algebra
relations, so
\begin{prop}[Relation (V)]
The relation from Figure~\ref{fig:Rel5} holds.
\end{prop}

Theorem~\ref{main1} from the Introduction follows.

\section{Functor valued invariants of oriented tangles}
\label{Section6}

\begin{figure}
  \centering
 \includegraphics{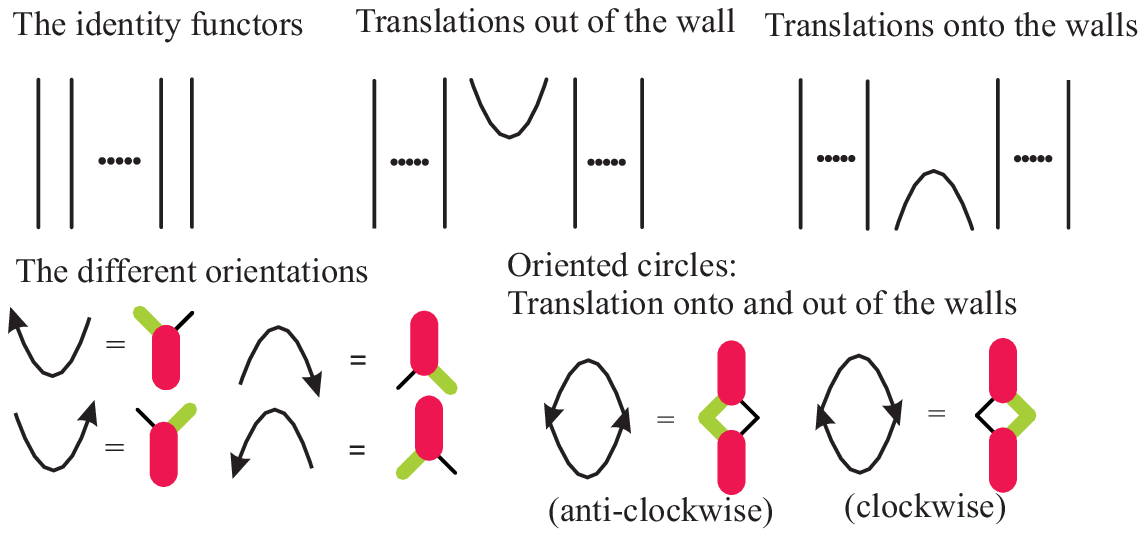}
  \caption{Crossingless elementary tangles and their associated functors}
  \label{tangles}
\end{figure}

We want to use the previous paragraphs to construct a functor valued
invariant of oriented tangles categorifying the quantum
$\mathfrak{sl}_k$-invariants.

If $\cA$ is an abelian category we denote by $\cD^b(\cA)$ the bounded derived
category with shift functor $\llbracket\;\rrbracket$ such that
$\llbracket1\rrbracket$ shifts the complex one step to the right.

Recall now the definition of the tangle category $\cT$ (see for example
~\cite{Kassel}, \cite{Tom}). The objects are finite ${+,-}$-sequences,
including the empty sequence; morphisms are the isotopy classes of oriented
tangles. Here a plus indicates the orientation downwards, whereas a minus
indicates the orientation upwards. The unoriented elementary tangles are
depicted at the top of Figure~\ref{tangles}. The first cup below would be a
morphism from the emptyset to $(-,+)$, whereas the cup in the left lower corner
is a morphism from the emptyset to $(+,-)$. Any morphism in $\cT$ is a
composition of oriented elementary morphisms.

For any object $a\in\cT$ we define $|a|:=j+(k-1)i$ where $i$ is the number of
pluses and $j$ the number of minuses in $a$. To an elementary morphism from $a$
to $b$ we associate a functor
$\cF:\cD^b(\bigoplus_\mu{}^\mZ\cO(|a|)^\mu)\rightarrow
\cD^b(\bigoplus_{\mu'}{}^\mZ\cO(|b|)^{\mu'})$, where $\mu$ and $\mu'$ run
through all partitions with at most $k$ parts, as follows:

\begin{figure}
  \centering
 \includegraphics{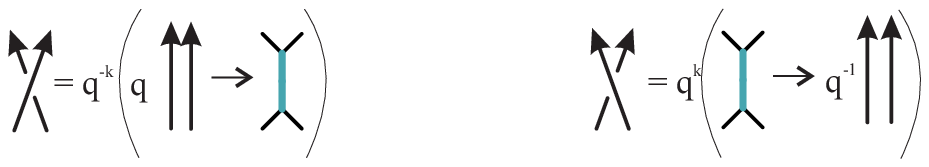}
  \caption{Crossingless elementary tangles and their associated functors}
  \label{crossings}
\end{figure}

\begin{enumerate}
\item To vertical strands we associate the identity functor
(Figure~\ref{tangles}) between the associated categories.
\item A cap diagram should first be replaced by a trivalent graph with labels $1$, $k-1$ and $k$,
depending on its orientation, and as shown in Figure~\ref{tangles}. To a cap
diagram we associate the corresponding standard lift of translation functor
onto the walls as defined in Section~\ref{translations}. The orientation
determines the corresponding categories (Figure~\ref{tangles}).
\item A cup diagram should first be replaced by a trivalent graph with labels $1$, $k-1$ and $k$,
depending on its orientation, and as shown in Figure~\ref{tangles}. To a cup
diagram we associate then the corresponding standard lift of translation
functor out of the walls as defined in Section~\ref{translations}.
\item Following \cite{StDuke}, we associate to a positive crossing with upwards pointing arrows the
corresponding left derived of the shuffling functor, but now shifted by
$\langle -k\rangle\llbracket1\rrbracket$. To a negative crossing we associate
the right derived of the coshuffling functor shifted by $\langle
k\rangle\llbracket-1\rrbracket$. In other words, we take the cone of the
natural transformations as depicted in Figure~\ref{crossings}, where the
identity parts are concentrated in position zero of the complex. The natural
transformations are both homogeneous of degree zero and arise as adjunction
morphisms from translation on and out of the wall.
\end{enumerate}

\begin{figure}
  \centering
 \includegraphics{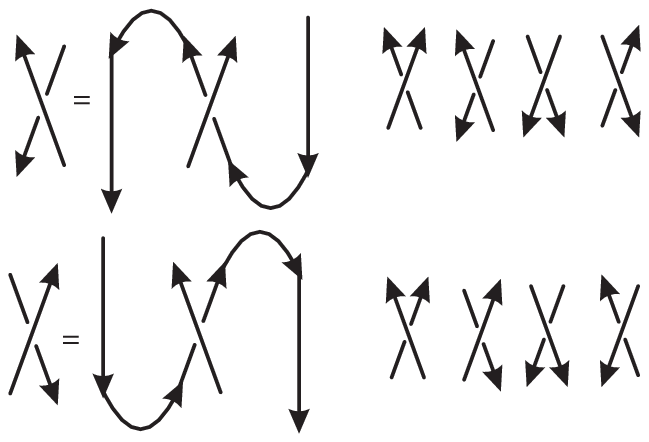}
  \caption{The functors associated with arbitrary crossings}
  \label{arbitrarycrossings}
\end{figure}

To an arbitrary crossing we associate the functors given in
Figure~\ref{arbitrarycrossings}: We first consider the positive upwards
pointing crossing and compose it with cap and cap as indicated to get the
negative crossing pointing to the left. Repeating this process we get all the 4
crossings depicted to the right in the first row of
Figure~\ref{arbitrarycrossings}. Analogously we could start with the (negative)
upwards pointing crossing and proceed as shown in the second row of
Figure~\ref{arbitrarycrossings}. This associates with each type of crossing a
functor. To obtain Theorem~\ref{main2} from the introduction we have to check
the invariance under tangle moves.

\subsection{The tangle moves} In Figure~\ref{Reidemeister1} we
have depicted four pairs of functors. In the first pair, the functor $F$ on the
RHS has been already defined and goes from the singularity $\nu=(1,k)$ to
$\nu=(k,1)$. The corresponding categories can be identified via an
Enright-Shelton equivalence (\cite{ES}). The following proposition ensures that
under this identification the functor $F$ becomes isomorphic to the identity
functor. We indicate the identifications to be made by slightly incline the
arrow. Analogous statements hold for the remaining three functors shown in
Figure~\ref{Reidemeister1}. Hence the following result should be considered as
a refined version of the isotopy relations of tangles:

\begin{figure}
  \centering
 \includegraphics{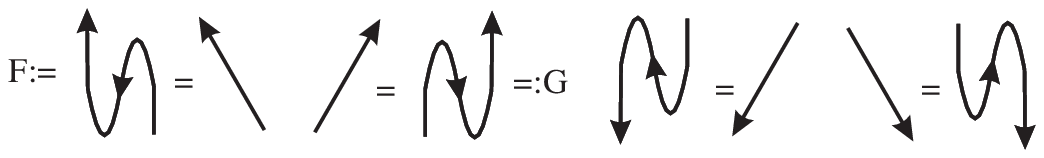}
  \caption{The 4 versions of the Isomorphism 1 of Tangles}
  \label{Reidemeister1}
\end{figure}

\begin{prop}[Isomorphisms 1]
\label{propReidemeister1}
The functors depicted in
Figure~\ref{Reidemeister1} are all equivalences of the corresponding
categories. (The first two functors are mutually inverse, as so are
the second two functors).
\end{prop}

\begin{proof}

Let $F'=\theta^{(1,k)}_{(1,k-1,1)}\theta_{(k,1)}^{(1,k-1,1)}$ and
$G'=\theta^{(k,1)}_{(1,k-1,1)}\theta_{(1,k)}^{(1,k-1,1)}$ be the naively
associated functors to the graphs of Figure~\ref{Reidemeister1}.
Combinatorially, the composition $G'F'$ is given as follows:
\begin{eqnarray*}
\begin{array}{l|l|l|l}
   (k,1)&(1,k-1,1)&(1,k)&(1,k-1,1)\\
\hline e&A_{k-1}^{1}&A_{k-1}^{1}&A_{k-1}^{1}B_{2}^{k}
\end{array}
\end{eqnarray*}
\begin{figure}
 \includegraphics{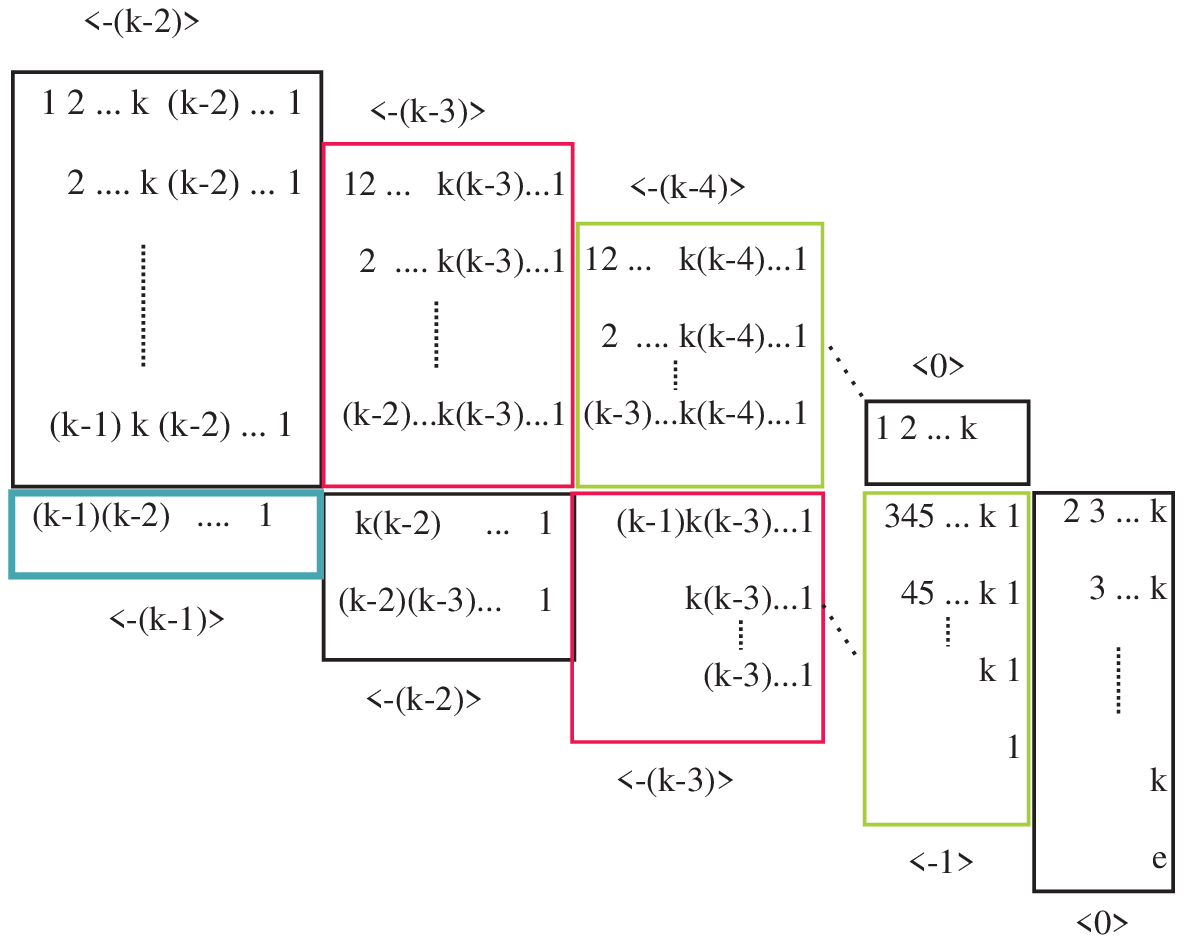}
 \caption{The Verma flag of $F'G'\Delta(e)$}
  \label{fig:boxes}
\end{figure}
The braid relations in $S_n$ provide the equality $$r\ldots 2\;1\;2\;3\;\ldots
r=1\;2\ldots r(r-1)\ldots 2\;1$$ for any $1\leq r< n.$ Using these equalities
one can show that $A_{k-1}^{1}B_{2}^{k}$ is of the form as depicted in
Figure~\ref{fig:boxes}. The top line of the $i$-th box upstairs is in degree
$i-1$, whereas the bottom line is always in degree $k-2$. The top line of each
box downstairs is in degree $k-1$, whereas the bottom line of the $i$-th box is
in degree $k-2+i$. We combine the $i$-th upstairs box with the $i+1$-th
downstairs box.

Translating to $(k,1)$, any two combined boxes together represent (up to a
shift in the grading) a copy of the projective module $P:=P(1\;2\;\ldots k)$.
(Above or below each box we denoted the grading shift $\langle -j\rangle$ which
occurs if we translate any element $x$ from the box to $(k,1)$ - one just has
to remove the last $j$ elements from $x$ and shift by $-j$ in the grading). The
only remaining element from the first downstairs box becomes a copy of $P(e)$.
Altogether we get $G'F'\Delta(e)=[k-1]P\oplus P(e)$. The projective module
$\theta_{k,1}^{(1^{k+1})}P$ corresponds to the following permutation (of $k+1$
letters)
\begin{eqnarray*}
x=\left(
\begin{array}[tbh]{cccccccccccccc}
1&2&3&\ldots&k-1&k&k+1\\
k+1&k&k-1&\ldots&3&2&1
\end{array}
\right)
\end{eqnarray*}

Under the Robinson-Schensted correspondence this corresponds to a tableau with
entries $1,2,\ldots, k+1$ in its first column. {\bf Fact 3} implies now
$FG\cong\operatorname{id}_{(k,1)}$. We leave it to the reader to verify that
$G'F'\Delta(e)\cong P(k\;(k-1)\ldots1)\oplus\Delta(e)$, where the first summand
translated out of the walls is $P(x)$, where $x$ is as above.
Invoking again {\bf Fact 3}, it follows $GF\cong\operatorname{id}_{(1,k)}$.
Hence the functors $F$ and $G$ define mutually inverse equivalences of (the
singular parabolic) categories in question. Similar calculations show that the
remaining two functors are mutually inverse equivalences as well, we omit the
details.
\end{proof}

\begin{figure}
  \centering
 \includegraphics{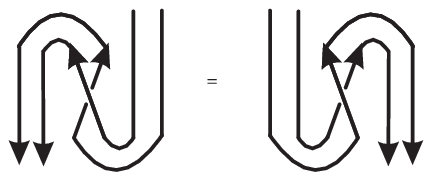}
  \caption{Isomorphism 2 of Tangles}
  \label{Reidemeister2}
\end{figure}

\begin{prop}
\label{PropReidemeister2} The functors associated to the tangle
diagrams depicted in Figure~\ref{Reidemeister2} are isomorphic.
\end{prop}

As preparation we need to prove several small statements, formulated
as Lemmas.

\begin{figure}
  \centering
 \includegraphics{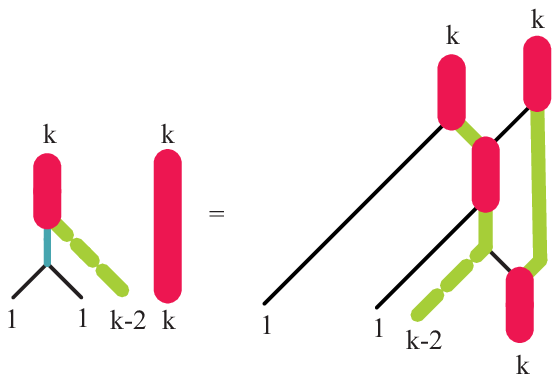}
 \caption{Step 1 in the proof of Proposition~\ref{PropReidemeister2}}
  \label{Lemma53}
\end{figure}

\begin{lemma}
There is an isomorphisms of functors as shown in
Figure~\ref{Lemma53}.
\end{lemma}

\begin{proof}
The proof is again completely combinatorial, so we leave out the details. The
functor associated with the left hand side maps $\Delta(e)$ to $\Delta(e)$. The
functor associated with the right hand side maps $\Delta(e)$ to a direct sum of
$\Delta(e)$ and copies of $P:=P\left((2k-1)\ldots(k+1)2\ldots k\right)$. On the
other hand $\theta_{(k,k)}^{(1^{2k})}P=P(x)$, where

\begin{eqnarray*}
x=\left(
\begin{array}[tbh]{ccccccccccccccc}
1&2&3&\ldots&k-1&k&k+1&k+2&\ldots&2k-1&2k\\
2k&k&k-1&\ldots&3&1&2k-1&2k-2&\ldots&k+1&2
\end{array}
\right),
\end{eqnarray*}
and so $x$  corresponds to a tableaux with the numbers $1,3,4,\ldots k,2k-1,2k$
in the first column, which means there are $k-1+2=k+1$ rows. The statement
follows by applying {\bf Fact 3}.
\end{proof}

\begin{figure}
  \centering
 \includegraphics{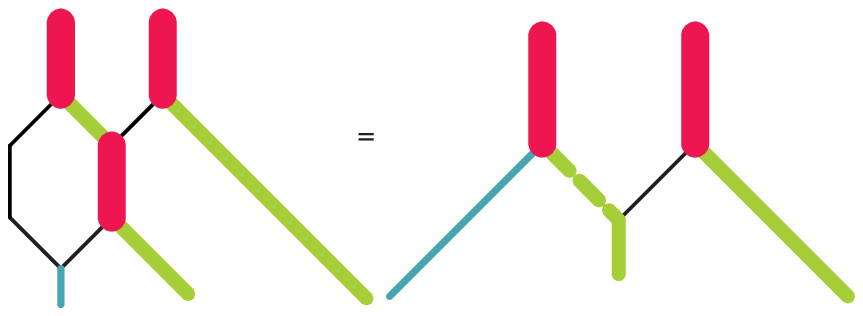}
 \caption{Step 2 in the proof of Proposition~\ref{PropReidemeister2}}
  \label{Lemma54}
\end{figure}

\begin{lemma}
\label{Step2} There is an isomorphisms of functors as shown in
Figure~\ref{Lemma54}.
\end{lemma}

\begin{proof}
The proof is again completely combinatorial, so we leave out the
details. The functor associated with the right hand side maps
$\Delta(e)$ to $P(3\ldots k)$, whereas the functor associated with
the left hand side maps $\Delta(e)$ to $P(12\ldots k)\oplus
P(3\ldots k)$. Note that $\theta_{(k,k)}^{1^{2k}}P(12\ldots
k)=P(x)$, where
\begin{eqnarray*}
x=\left(
\begin{array}[tbh]{cccccccccccccc}
1&2&\ldots&k&k+1&\ldots&2k-1&2k\\
k+1&k&\ldots&2&2k&\ldots&k+2&1
\end{array}
\right),
\end{eqnarray*}
and so $x$  corresponds to a tableaux with the numbers
$1,2,\ldots,k+1$ in the first column, which means there are $k+1$
rows. The statement follows.
\end{proof}

\begin{proof}[Proof of Proposition~\ref{PropReidemeister2}]
Let $F_1$ (resp. $F_2$) be the functor on the left (right) hand side
of Figure~\ref{Lemma53}. Let $G_1$ (resp. $G_2$) be the functor on
the left (right) hand side of Figure~\ref{Lemma54}. Fix any
composition $\mu$ of $2k$ with at most $k$ parts and consider the
functors
\begin{eqnarray*}
\begin{array}[thb]{lllllllll}
H&:=&\theta_{(1,1,k-1,k-1)}^{(2,k-1,k-1)}:&&\cO_{(1,1,k-1,k-1)}^\mu&\rightarrow&\cO_{(2,k-1,k-1)}^\mu\\
J&:=&\theta^{(1,1,k-2,k)}_{(1,1,k-2,1,k-1)}\theta_{(1,1,k-1,k-1)}^{(1,1,k-2,1,k-1)}:
&&\cO_{(1,1,k-1,k-1)}^\mu&\rightarrow&\cO_{(1,1,k-2,k)}^\mu
\end{array}
\end{eqnarray*}

Then we have isomorphisms of functors as follows:
\begin{eqnarray*}
G_1H\cong G_2H\cong F_1J\cong F_2J.
\end{eqnarray*}
This follows directly from the Lemmas~\ref{Lemma53} and~\ref{Lemma54} by
drawing pictures. Using Proposition~\ref{Reidemeister1},
Figure~\ref{Reidemeister1} we see that $F_2J$ is isomorphic to the functor
given by the vertically reflected diagram. From this it follows that we have an
isomorphism of functors as in Figure~\ref{Reidemeister2}, but the crossings
replaced by \includegraphics{theta.eps}. Now one has just to take the Cone of
the corresponding adjunction morphism. Up to a scalar, there is a unique
morphism with the correct degree. The statement of the
Proposition follows by applying {\bf Fact 3}.
\end{proof}

\begin{prop}[Reidemeister 2 and 3]
The functors associated to the positive and negative upwards
pointing crossings are mutually inverse equivalences and satisfy the
braid relations.
\end{prop}

\begin{proof}
This is a standard fact, see for example \cite{MS}.
\end{proof}

\begin{figure}
\centering
 \includegraphics{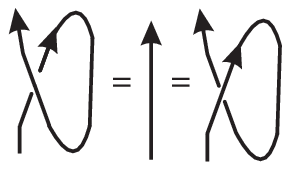}
  \caption{Reidemeister 1}
  \label{Reidemeister5}
\end{figure}

\begin{prop}[Reidemeister1]
\label{PropRel5} The three functors associated to the tangle diagrams in
Figure~\ref{Reidemeister5} are isomorphic.
\end{prop}

\begin{proof}
The functors in question are going from the singularity $(1,k)$ to
the singularity $(1,k)$. Recall the definition of the functor
associated to the crossings.

Let us first give a short explanation why one might expect the
claimed isomorphisms: From the relations in Figures~\ref{fig:Rel3}
and Figure~\ref{fig:Rel2} the functor on the left hand side of
Figure~\ref{Reidemeister5} is,
up to an overall shift by $\langle -k\rangle$, the Cone of a
morphism
$$\gamma: \quad [k]\operatorname{id}\langle 1\rangle  \longrightarrow [k-1]\operatorname{id},$$
sitting in cohomological degree zero and $1$. There is
the obvious surjection
$$(q^{k}+q^{k-2}+q^{k-4}+\ldots+q^{4-k}+q^{2-k})\operatorname{id}\rightarrow
(q^{k-2}+q^{k-4}+\ldots+q^{4-k}+q^{2-k})\operatorname{id}$$ which identifies
the same summands and has kernel $q^{k}\operatorname{id}$, so that we expect
the second isomorphism of Figure~\ref{Reidemeister5} (and similarly the first
one). To prove the statement we have to understand the morphism $\gamma$
better.

The adjunction morphism
$\alpha:\theta_{(1,k)}^{(1,1,k-1)}\rightarrow\theta_{2,k-1}^{1,1,k-1}\theta_{1,1,k-1}^{2,k-1}
\theta_{(1,k)}^{(1,1,k-1)}$ is injective for any module with Verma
flag, in particular for Verma modules and projectives. From the
proof of Proposition~\ref{PropRel3} we see that the image of the
adjunction morphism applied to $\Delta(e)$ is a module with Verma
subquotients given by $qA_2^k$, $q^k e$. Hence
$\gamma':=\theta_{(1,1,k-1)}^{(1,k)}(\alpha)_{\Delta(e)}$ surjects
onto the $[k-1]$ copies of $\Delta(e)$, and defines a split
$$\theta^{(1,k)}_{(1,1,k-1)}\theta_{(1,k)}^{(1,1,k-1)}\cong F'\oplus
[k-1]\operatorname{id}$$ for some projective functor $F'$. Thanks to
Proposition~\ref{propRel2} we have $F'\cong\operatorname{id}\langle
k\rangle$. Now, if we restrict to the parabolic subcategories with
at most $k$ parts, then $\gamma'$ induces the surjection with kernel
the identity functor shifted up by $k$ in the degree. Putting the
overall shift back into the picture, we obtain the second
isomorphism. The first isomorphism can be proved analogously or by
observing that these are just the adjoint functors.
\end{proof}

\begin{figure}
\centering
 \includegraphics{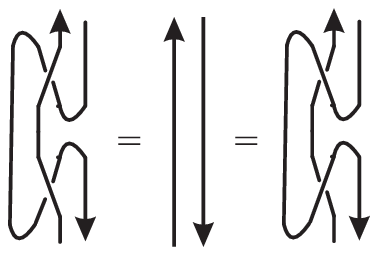}
 \hspace{1.5cm}
 \includegraphics{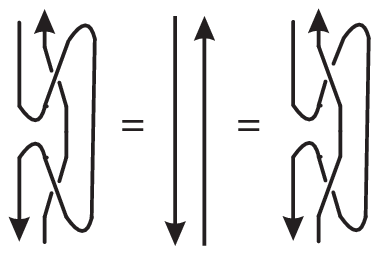}
  \caption{Isomorphisms of Tangles}
  \label{Reidemeister6}
\end{figure}

\begin{prop}
\label{Sussangap} The functors associated to the tangle diagrams in
Figure~\ref{Reidemeister6} satisfy the displayed isomorphisms.
\end{prop}

\begin{proof}
The right half of Figure~\ref{Reidemeister6} is just the reflection in a
vertical line of the diagrams in Figure~\ref{Reidemeister6}. Now there is an
isomorphism of the Lie algebra $\mathfrak{gl}_n$ given by the obvious
involution of the Dynkin diagram which swaps the $i$-th with the $n-1-i$-th
node. This isomorphism defines an auto-equivalence of the category $\cO$ for
$\mathfrak{gl}_n$ which identifies ${}^\mZ\cO(n)_\nu^\mu$ with
${}^\mZ\cO_{\tilde\nu}^{\tilde{\mu}}$, where the partition are ``reflected in a
vertical line''. Under this automorphism the functors displayed on the left half
of Figure~\ref{Reidemeister6} correspond to the functors displayed on the right
half, so that it is enough to prove the first two isomorphisms.\footnote{Note
that the proof of the corresponding result in \cite{Sussan} is not complete.}
Consider first the diagram on the left hand side together with the following
functors
\begin{eqnarray*}
F:=\theta_{(k,1,k-1)}^{(k-1,1,1,k-1)},&&\hat{F}:=\theta^{(k,1,k-1)}_{(k-1,1,1,k-1)}\\
H:=\theta_{(k-1,1,1,k-1)}^{(k-1,1,k)},&& \hat{H}:=\theta^{(k-1,1,1,k-1)}_{(k-1,1,k)},\\
 &&\theta:=\theta^{(k-1,1,1,k-1)}_{(k-1,2,k-1)}\theta_{(k-1,1,1,k-1)}^{(k-1,2,k-1)}\\
 G:=\operatorname{Cone}(\theta\longrightarrow\operatorname{id}\langle
 -1\rangle)\langle k\rangle,&&  \hat{G}:=\operatorname{Cone}(\operatorname{id}\langle
 1\rangle\longrightarrow\theta)\llbracket1\rrbracket\langle -k\rangle,
\end{eqnarray*}
The relation we want to verify says exactly that after restricting to
parabolics with at most $k$ parts, the functors
$\Phi_1:=\hat{F}\,\hat{G}\,\hat{H}$ and $\Phi_2:=H\,G\,F$ are inverse to each
other.

Directly from the definitions it follows, that the composition $\Phi_1\Phi_2$
is given by the the following complex of functors:
\begin{eqnarray*}
\hat{F}\hat{H}\,H\,\theta\,F\langle 1\rangle\longrightarrow
\hat{F}\,\theta\,\hat{H}\,H\,\theta\,F \oplus \hat{F}\hat{H}\,H\,F
\longrightarrow \hat{F}\,\theta\,\hat{H}\,H\,F\langle-1\rangle.
\end{eqnarray*}
Here the first map is $\left(\begin{array}{c}
\hat{F}(\beta)_{\hat{H}H\theta F}\\
\hat{F}\hat{H}H(\alpha)_{F}\end{array}\right)$, and the second is
$\big(\hat{F}\theta\hat{H}H(\alpha)_{F},-\hat{F}(\beta)_{\hat{H}HF}\big)$,
where $\alpha$
is the adjunction morphism $\theta\rightarrow\operatorname{id}\langle-1\rangle$
and $\beta$ the adjunction morphism $\operatorname{id}\langle
1\rangle\rightarrow\theta$.

Using now the Relations (I), (III) and (IV) (Figures~\ref{fig:Rel1},
\ref{fig:Rel3}, \ref{fig:Rel4}) the restrictions of the functors to any
parabolic with at most $k$ parts gives rise to the complex
\begin{eqnarray}
\label{complex}
[k-1]J\langle1\rangle\longrightarrow(\operatorname{id}\oplus[k-2]J)\oplus
[k]J\longrightarrow[k-1]J\langle-1\rangle
\end{eqnarray}
where $J$ is the restriction of the functor
$\theta^{(k,k-1,1)}_{(k,k)}\theta_{(k,k-1,1)}^{(k,k)}$. As in
Proposition~\ref{PropRel5} we deduce that the first map is an inclusion
and the second map is a surjection so that the functor
$[k-1]J\langle-1\rangle$ splits off as a direct summand and \eqref{complex} is
quasi-isomorphic to
\begin{eqnarray}
\label{complex2}
0\rightarrow [k-1]J\langle1\rangle
\overset{\gamma}{\longrightarrow}\operatorname{id}\oplus[k-1]J\rightarrow 0.
\end{eqnarray}
Denote by $\kappa_1:\operatorname{id}\oplus[k-1]J\rightarrow [k-1]J$ the
projection. We claim that $\kappa_1\gamma$ is an isomorphism.

Indeed, assume that this is not the case. Let $P(w)$ be
an indecomposable projective, different from the
dominant Verma module $\Delta(e)$.
Then the Verma flag of $P(w)$ contains, as a submodule, the
copy of $\Delta(e)$ which corresponds to the
inclusion $\Delta(w)\hookrightarrow\Delta(e)$. The
socle of this submodule is in the kernel of any non-invertible
homomorphism $f:P(w)\rightarrow P(w)$ and
any homomorphism $g:P(w)\rightarrow\Delta(e)$. Thus
it is in the kernel of $\gamma_{\Delta(e)}$, which contradics the
injectivity of $\gamma$.

Let $\kappa_2:\operatorname{id}\oplus[k-1]J\rightarrow
\operatorname{id}$ be the projection. In particular,
$\gamma=\left(\begin{array}{c}\kappa_2\gamma\\
\kappa_1\gamma\end{array}\right)$. Now the map
$(\mathrm{id},-\kappa_2\gamma(\kappa_1\gamma)^{-1}):
\operatorname{id}\oplus[k-1]J\rightarrow \operatorname{id}$ satisfies
$(\mathrm{id},-\kappa_2\gamma(\kappa_1\gamma)^{-1})\gamma=0$ and hence defines
a quasi-isomorphism from \eqref{complex2} to the complex
$0\rightarrow\operatorname{id}\rightarrow 0$, which represents the identity
functor. This proves the first isomorphism. The second can be deduced
analogously. Alternatively one could deduce it by adjointness properties.
\end{proof}

To summarise: Theorem~\ref{main2} from the introduction holds.

\section{Cohomology rings, natural transformations and foams}

In this final section we indicate how to extend our functorial
invariant of trivalent graphs to an invariant of trivalent graphs
and foams, and also explain the connection with \cite{Ksl3}.
Conjecturally our setup actually gives the representation
theoretical background for the very recent generalisation
\cite{MackaayStosic} of \cite{Ksl3} to arbitrary $k$.

Roughly speaking, {\it a foam} is a morphism between certain trivalent graphs
(for a precise definition see \cite{Ksl3}, \cite{MackaayVaz},
\cite{MackaayStosic}). Khovanov associated to each special trivalent graph a
graded vector space and to any foam a homogeneous linear map of degree being
the degree of the foam. In the following we want to indicate how this
construction emerges naturally from our picture by restricting the functors to
the non-parabolic part and applying some Soergel's combinatorial functor $\mV$.
In the following we assume that the reader is familiar with \cite{Ksl3}.

\subsection{Natural transformation associated with basic foams}
\begin{figure}
  \centering
 \includegraphics{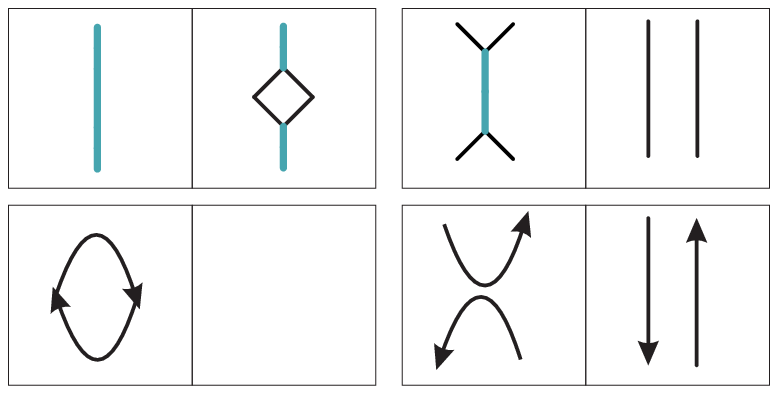}
  \caption{Basic Foams correspond to basic natural transformations, homogeneous of degree $-1$, $1$, $-2$, $2$ respectively.}
  \label{basicfoams}
\end{figure}

Apart from the identity morphisms, un-dotted foams are compositions of
elementary foams as depicted in Figure~\ref{basicfoams}.
Each rectangle should be read from the left to the right, as well as from the
right to the left; giving rise to two basic foams. Additionally, both possible
orientation should be considered in the last two cases. For each graph
appearing as the boundary of a foam, we have the associated functor
(Section~\ref{Section5}). We assign now to each basic foam a natural
transformation, all of them will be just adjunction morphisms:

{\it First row:} We associate the adjunction morphism $\beta_1$ from the
identity to the composition $\theta_{(1,1)}^{(2)}\theta_{(2)}^{(1,1)}$, and
$\beta_2$ vice versa. $\beta_1$ and $\beta_2$ are homogeneous of degree $-1$
(\cite[Theorem 8.4]{Stgrad}). Thanks to (\cite[Remarks 3.8 c)]{Stgrad})
we have adjunction morphisms
$\alpha_1:\operatorname{id}\rightarrow\theta_{(2)}^{(1,1)}\theta_{(1,1)}^{(2)}$
and
$\alpha_2:\theta_{(2)}^{(1,1)}\theta_{(1,1)}^{(2)}\rightarrow\operatorname{id}$,
both homogeneous of degree $1$. A priori, they are unique up to a non-zero
scalar - which we want to choose such that Lemma~\ref{adj} and
Lemma~\ref{adj2} below hold; the
same will apply to all the other adjunction morphisms. These are the natural
transformations we associate to the two foams given by the first diagram.

{\it The second row:} Recall that we associated to a circle the composition of
translation out of the walls and onto the walls as depicted in
Figure~\ref{tangles}. Hence we have the obvious adjunction morphisms $\gamma_1$
from a clockwise circle, $\gamma_2$ from an anticlockwise circle, $\gamma_3$ to
a clockwise circle, $\gamma_4$ to an anticlockwise circle. They are all
homogeneous of degree $1-k$. This follows from the adjunction
$(\theta^{(i,j)}_{(k)},\theta^{(k)}_{(i,j)})$, where
$(i,j)\in\{(1,k-1),(k-1,1)\langle 1-k\rangle\}$ (a special case of
\cite[Proposition 4.2]{FKS}). The adjunction morphisms
$\delta_1:\theta_{(3)}^{(2,1)}\theta_{(2,1)}^{(3)}\rightarrow\operatorname{id}$,
$\delta_1':\theta_{(3)}^{(1,2)}\theta_{(1,2)}^{(3)}\rightarrow\operatorname{id}$
and
$\delta_2:\operatorname{id}\rightarrow\operatorname\theta_{(3)}^{(2,1)}\theta_{(2,1)}^{(3)}$,
$\delta_2':\operatorname{id}\rightarrow\operatorname\theta_{(3)}^{(1,2)}\theta_{(1,2)}^{(3)}$,
are homogeneous of degree $k-1$ (by the combinatorics of
Section~\ref{Section4}).

From now on we stick to the case $k=3$ and illustrate the connection to
\cite{Ksl3}. Denote by $\deg{\mathtt{F}}$ the degree of a basic foam
$\mathtt{F}$. From the
definition it follows:

\begin{lemma}
Let $k=3$. For a basic foam $\mathtt{F}$ let $\phi_{\mathtt{F}}$
be the associated
natural transformation as defined above.
Then $\operatorname{deg}(\mathtt{F})=
\operatorname{deg}(\phi_{\mathtt{F}}).$
\end{lemma}

\begin{figure}
\centering
 \includegraphics{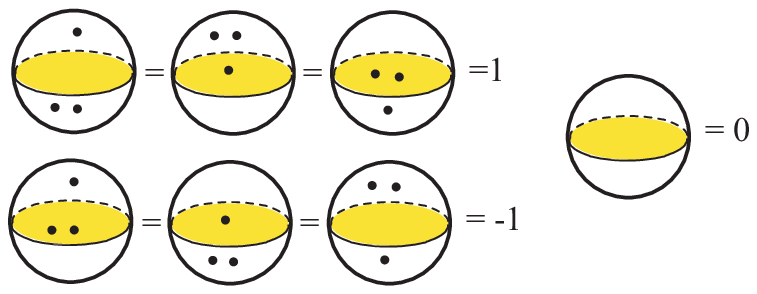}
  \caption{Examples of Theta foams with all the non-trivial evaluations}
  \label{thetafoams}
\end{figure}

Apart from the basic foams we need the so-called theta foams. Theta-foams
(Figure~\ref{thetafoams}) are obtained by gluing three oriented disks along
their boundaries (their orientations must coincide).

Dots will correspond to multiplication with a certain element of degree two in
the centre of the category. This will be exactly as in \cite{Ksl3} and
\cite{MackaayStosic}. To explain this connection we have to bring cohomology
rings of partial flag varieties into the picture.

\subsection{The cohomology of flag varieties}

Recall the following result of Soergel: The category ${}^\mZ\cO(n)_\la$
(for $\la$ a partition of $n$) has one indecomposable projective-injective
module $P_\la$ with head concentrated in degree zero. We have Soergel's
functor
$$\mV=\HOM_{\cO}(P(\la), _-):{}^\mZ\cO(n)_\la
\longrightarrow \gMOD-\END_{\cO}(P_\la).$$
By Soergel's Endomorphismensatz (\cite{Sperv}) we know that
$\END_{\cO}(P_\la)$ is isomorphic (as a graded ring) to the
cohomology ring (with complex
coefficients) of the associated partial flag variety $\cF_\la$, where the
dimensions of the subquotients are $\{\la_i\}_{i>0}$.

For instance $\END(P_{(1,1)})\cong H^*(\cF_{(1,1)})\cong\mC[x]/(x^2)$. If we
choose $\la=(3)$, then we just get the cohomology $\mC$ of a point, whereas
$\END_{\cO}(P_\la)\cong H^*(\mP^2)\cong\mC[x]/(x^3)$
if $\la=(2,1)$ or $\la=(1,2)$.
In each case, $x$ is of degree two. If we choose the reversed standard
orientation on $\mP^2$, then the cohomology ring $\cA:=\mC[x]/(x^3)$ comes
along (\cite{Ksl3}) with the trace form $\operatorname{Tr}(x^i)=-1\delta_{2,i}$
and the comultiplication $$\Delta(1)=-(1\otimes x^2+x\otimes x+x^2\otimes
1),\quad\Delta(x)=-(x\otimes x^2+x^2\otimes x),\quad \Delta(x^2)=-x^2\otimes
x^2.$$ We choose the basis $X_{(1)}=1,X_{(2)}=x,X_{(3)}=x^2$ of $\cA$ and
denote by $X^{(1)}$, $X^{(2)}$, $X^{(3)}$ its dual basis with respect to
$\operatorname{Tr}$.

Finally, the cohomology ring $C:=H^*(\cF_{(1,1,1)})$ is isomorphic to the
polynomial ring $\mC[X_1,X_2,X_3]$ modulo the ideal generated by the elementary
symmetric polynomials. There is the trace function $\operatorname{Tr}:
C\rightarrow \mC$ which maps $X_1X_2^2$ to $1$.

\subsection{The bridge}
The functor $\mV$ connects category $\cO$ and modules over cohomology rings of
flag varieties: The functor $\theta_{(2)}^{(1,1)}\theta^{(2)}_{(1,1)}:
{}^\mZ\cO(2)_{(1,1)}\rightarrow{}^\mZ\cO(2)_{(1,1)}$ corresponds (\cite{Sperv},
\cite{Stgrad}) under $\mV$ to the functor
$$\bullet\otimes_\mC\mC[x]/(x^2)\langle -1\rangle:\quad\gMOD-\mC[x]/(x^2)\rightarrow \gMOD-\mC[x]/(x^2).$$

\begin{lemma}
\label{adj} Under the above correspondence the natural
transformations $\alpha_1$, $\alpha_2$ become the
multiplication $\mV(\alpha_1)_N:
(N\otimes_\mC\mC[x]/(x^2)\langle -1\rangle)\langle 1\rangle
\rightarrow N$, $n\otimes c\mapsto nc$ and the comultiplication
$\mV(\alpha_2)_N: N\langle 1\rangle\rightarrow
N\otimes_\mC\mC[x]/(x^2)\langle -1\rangle$,
$n\mapsto x\otimes n+1\otimes xn$ respectively.
\end{lemma}
\begin{proof}
See \cite[Lemma 8.2]{StDuke}.
\end{proof}

Similarly, the functor $\theta^{(2)}_{(1,1)}\theta_{(2)}^{(1,1)}:
{}^\mZ\cO(2)_{(2)}\rightarrow {}^\mZ\cO(2)_{(2)}$ corresponds under $\mV$ to
the functor
\begin{eqnarray}
\label{outon} \bullet\otimes_\mC\mC[x]/(x^2)\langle
-1\rangle\cong\operatorname{id}\langle1\rangle\oplus
\operatorname{id}\langle-1\rangle:\quad
\gMOD-\mC\rightarrow\gMOD-\mC.
\end{eqnarray}

\begin{lemma}\label{adj2}
With the above definitions,
for every graded $\mC$-module $N$ we have $\mV(\beta_1)_N:
N\langle -1\rangle
\rightarrow N\otimes_\mC\mC[x]/(x^2)\langle -1\rangle,$ $ n\mapsto
n\otimes 1$ and we further  have the following: $\mV(\beta_2)_N:(N\otimes_\mC\mC[x]/(x^2)\langle
-1\rangle)\langle -1\rangle\rightarrow N,$ $ n\otimes c\mapsto
\operatorname{Tr}(c)n$. Under the isomorphism \eqref{outon} we just get the
projection and inclusion morphisms of degree $-1$.
\end{lemma}

\begin{proof}
Since the source and target categories of the functors are semi-simple, there
is only one (up to scalar) possible map of the correct degree in each case.
\end{proof}

Now consider the functor $\theta_{(3)}^{(2,1)}\theta^{(3)}_{(2,1)}:
{}^\mZ\cO(3)_{(2,1)}\rightarrow{}^\mZ\cO(3)_{(2,1)}.$ Under the functor $\mV$
this corresponds to the functor
\begin{eqnarray}
\label{eq} \bullet\otimes_\mC\mC[x]/(x^3)\langle
-2\rangle&\cong&\bullet\otimes_\mC\cA\langle -2\rangle
\end{eqnarray}
(\cite[3.4]{FKS}). Because of Soergel's double centralizer property with
respect to the antidominand projective module, a natural transformation between
projective functors is already determined by its value on the antidominant
projective module (by argumens similar to e.g.
\cite[Lemma~5.1]{MazStr}). Hence the following Lemma is useful:

\begin{lemma}
Under the functor $\mV$ we have the following correspondences:
\begin{itemize}
\item Evaluated at the antidominant projective module $P(21)$ or $P(12)$, the
natural transformations $\delta_1$ and $\delta_1'$ correspond to the
multiplication morphism $m:(\cA\otimes\cA\langle -2\rangle)
\langle 2\rangle\rightarrow\cA$,
whereas $\delta_2$ and $\delta_2'$ corresponds to the comultiplication
morphism $\Delta:\cA\langle 2\rangle\rightarrow\cA\otimes\cA\langle -2\rangle$.
\item Evaluated at the dominant Verma
module $\Delta(e)$, we get for $\delta_1$ and $\delta_1'$ the induced
multiplication morphism $\overline{m}:(\mC\otimes\cA
\langle -2\rangle)\langle 2\rangle\rightarrow\mC$,
and for $\delta_1$ and $\delta_1'$ the induced comultiplication morphism
$\overline\Delta:\mC\langle 2\rangle\rightarrow\mC\otimes\cA\langle -2\rangle$.
\end{itemize}
\end{lemma}

\begin{proof}
Note first that we have
$\mV\theta_{(3)}^{(2,1)}\theta_{(2,1)}^{(3)}P(21)\cong\cA\otimes\cA\langle-2\rangle$,
and similarly
$\mV\theta_{(3)}^{(1,2)}\theta_{(1,2)}^{(3)}P(12)\cong\cA\otimes\cA\langle
-2\rangle$ by \eqref{eq}; whereas
$\mV\theta_{(3)}^{(2,1)}\theta_{(2,1)}^{(3)}\Delta(e)\cong\mC\otimes\cA\langle -2\rangle$, and
similarly
$\mV\theta_{(3)}^{(1,2)}\theta_{(1,2)}^{(3)}\Delta(e)\cong\mC\otimes\cA\langle -2\rangle$.
Frobenius reciprocity provides a natural isomorphism of the form
\begin{eqnarray*}
\HOM_{\gMOD-\cA}(N\otimes\cA,N)\cong\HOM_{\gMOD-\mC}(N,N)
\end{eqnarray*}
mapping $f$ to $\hat{f}$, where $\hat{f}(n)=f(1\otimes n)$ for any graded right
$\cA$-module $N$ and $n\in N$. In particular, $\hat{m}$ is the identity map
which implies half of the statement.

Denote by $X^*$ the graded vector space dual of $X$. Then there is an
isomorphism of graded right $\cA$-modules as follows:
\begin{eqnarray*}
\gamma:\quad (N\otimes\cA)^*&\cong&N^*\otimes \cA, \quad
g\mapsto\sum_{i=1}^3g_i\otimes X^{(i)},\quad \widetilde{f\otimes c}
\mapsfrom f\otimes c,
\end{eqnarray*}
where $g_i(n)=g(n\otimes X_{(i)})$ and $\widetilde{f\otimes c}(n\otimes
d)=\operatorname{Tr}(cd)f(n)$ for $n\in N$, $c,d\in\cA$. The second adjunction
morphism is then the chain of isomorphisms
\begin{eqnarray*}
\HOM_{\gMOD-\mC}(N,N)&\cong&\HOM_{\gMOD-\mC}(N^*,N^*)\cong
\HOM_{\gMOD-\cA}(N^*\otimes\cA,N^*)\\
&\cong&\HOM_{\gMOD-\cA}((N\otimes\cA)^*,N^*)\cong
\HOM_{\gMOD-\cA}(N,N\otimes\cA).
\end{eqnarray*}
The first isomorphism here is the duality, the second the adjunction from
above, then we invoke the isomorphism $\gamma$ and finally the duality again.
It is now an easy direct calculation to verify the claim.
\end{proof}

\begin{lemma}
Under the functor  $\mathbb{V}$ for every graded $\mC$-module $N$ we have
the following:
$\mV(\gamma_1)_N: N\langle
-2\rangle\rightarrow N\otimes_\mC\mC[x]/(x^3)\langle -2\rangle, n\mapsto
n\otimes 1$ and further we have
$\mV(\gamma_2)_N:(N\otimes_\mC\mC[x]/(x^3)\langle
-2\rangle)\langle -2\rangle\rightarrow N, n\otimes c\mapsto
\operatorname{Tr}(c)n$. Under the isomorphism from Figure~\ref{fig:Rel1} we
just get the inclusion and projection morphisms of degree $-2$. The same
holds for $\gamma_1'$ and $\gamma_2'$.
\end{lemma}

\subsubsection{Dots on basic foams}
We still have to explain what to do with dots on basic foams. Under the functor
$\mV$ any dot just corresponds to multiplication with the variable $x$. By
Soergel's Struktursatz this means that we multiply the natural transformation
with a certain element of the centre of one of the involved categories
(\cite{Sperv}, \cite{MazStr}). To make this explicit, consider the functors
$F:=\theta_{(3)}^{(2,1)}\theta^{(3)}_{(2,1)}$ and
$G:=\theta_{(3)}^{(1,2)}\theta^{(3)}_{(1,2)}$. A natural transformation $f:
F\rightarrow F$ (or $g:G\rightarrow G$) is uniquely determined by $\mV(f):\mV
F\Delta(e)\rightarrow \mV F\Delta(e)$ (or $\mV(g):\mV G\Delta(e)\rightarrow \mV
G\Delta(e)$) (because $\cO(3)_{(3)}$ is semisimple).

\begin{figure}[htb]
\begin{center}
 \includegraphics{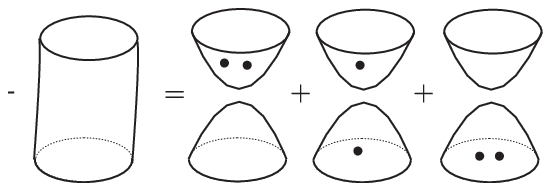}
\end{center}
  \caption{The surgery relation decomposes the identity morphisms}
  \label{fig:surgery}
\end{figure}

Choosing for $f$ and $g$ the identity morphism, we have $\mV(f), \mV(g):
\cA\rightarrow\cA$, and one checks directly that the surgery operation from
Figure~\ref{fig:surgery} decomposes them as follows:
\begin{eqnarray*}
-\operatorname{id}&=&m_x\circ
m_x\circ\delta_2\circ\delta_1+m_x\circ\delta_2\circ\delta_1\circ
m_x+\delta_2\circ\delta_1\circ m_x\circ m_x\\
-\operatorname{id}&=&m_x\circ
m_x\circ\delta_2'\circ\delta_1'+m_x\circ\delta_2'\circ\delta_1'\circ
m_x+\delta_2'\circ\delta_1'\circ m_x\circ m_x
\end{eqnarray*}
where $m_x$ is the multiplication with $x$ which we associate with a dot.

\subsubsection{Theta foams}
We have to associate to each theta foam a natural transformation from the
identity functor on ${}^\mZ\cO(3)_{(3)}$ to itself. To a theta foam with $d_i$
dots on the $i$-th disk we associate the natural transformation which
corresponds under the functor $\mV$ to the map $\mC\rightarrow\mC$,
$z\mapsto\operatorname{Tr}(X_1^{d_1}X_2^{d_2}X_3^{d_3})z$. In particular,
corresponding to the three discs (the equatorial, the upper hemisphere and the
lower hemisphere) there are three embeddings of $\cA$ into $C$, namely
$x\mapsto X_1$, $x\mapsto X_2$ and $x\mapsto X_3$ and we apply the usual rule
for the dots.

Let $\mathtt{F}$ be a basic foam with input boundary
$D_{F_1}$ and output boundary
$D_{F_2}$. Let $F_1$, $F_2$  be the corresponding functors as assigned in
Section~\ref{Section4} and  $G_1$, $G_2$ the associated graded vector spaces in
\cite{Ksl3}. Assigned to $\mathtt{F}$ we have
$\phi_{\mathtt{F}}:F_1\rightarrow F_2$
and also a linear map $g:G_1\rightarrow G_2$ from \cite{Ksl3}. Let
$\overline{F}_1$, $\overline{F}_2$ and $\overline{\phi_{\mathtt{F}}}$
be the restrictions to the non-parabolic summand. The following
result is now easily verified:

\begin{prop}
\begin{enumerate}
\item The above assignments define a functor from the category of prefoams as
defined in \cite{Ksl3} to the category of graded projective functors associated
with intertwiners and natural transformations between them.
\item There are isomorphism
$\mV \overline{F}_i\Delta(e)\cong G_i$, $i=1,2$, of graded vector spaces under
which $\mV\overline{\phi_{\mathtt{F}}}$ corresponds to $g$.
\end{enumerate}
\end{prop}

In particular, the approach of \cite{Ksl3} follows directly from our setup by
restriction. Note that we really loose some information here, since we evaluate
the natural transformation on the dominant Verma module (instead of on the
antidominant projective which would keep all the information). On the other
hand, we restricted to a direct summand. This is irrelevant for the quality of
the invariant, but only carries the information of the zero weight space in our
original $\mathfrak{sl}(k)$-modules $X^\nu$.

\begin{conjecture}
\label{conjecture} {\rm The obvious generalisation of our construction for
general $k$ gives rise to the Mackaay-Stosic-Vaz homology
(\cite{MackaayStosic}) and hence to the Khovanov-Rozansky homology
\cite{KhovRozMatrix}}.
\end{conjecture}

A verification of this conjecture would in particular imply a very nice
description of the interplay of natural transformations between projective
functors in terms of Schur polynomials, based on \cite{MackaayStosic}.

\subsection{Speculations on web bases and dual canonical bases}
In Section~\ref{Section4} we associated to each special intertwiner or
web diagram a certain projective functor. In the case $k=2$ the web
bases coincides with the Temperley-Lieb algebra basis which agrees with
Lusztig's canonical basis (\cite{FK}). One can show that the associated
functors are all indecomposable (\cite{StDuke}). This is however just
pure accident and very special for $k=2$. The answer to the following
question might shed some light on the relationship in general:

{\bf Question.} Is it true that the transformation matrix between the
web basis and the canonical basis describes the decomposition of the
functors assigned to webs into a direct sum of indecomposable functors?

To answer this question one has to improve the categorification
presented in the present paper, to include more general intertwiners,
and then connect it with the results on dual canonical bases from
\cite{FKK}, and the more general results \cite{BrundandualcanKL}.
Since there is no classification of indecomposable projective functors
for parabolic categories we expect that finding an answer to this
question might be quite hard.


\bibliography{ref}
\bibliographystyle{plain}
\end{document}